\documentclass[letter,twocolumn]{autart}    
                                     
\usepackage{graphics} 
\usepackage{epsfig} 
\usepackage{amsmath} 
\usepackage{amssymb}  
\usepackage{color}
\usepackage{natbib}
\usepackage{enumerate}
\usepackage{graphicx}
\usepackage{mathtools}
\usepackage{cancel}
\usepackage{url}

\usepackage[dvipsnames, table]{xcolor}
\usepackage{bbm}
\usepackage{tikz}
\usepackage{pgfplots}
\usepgfplotslibrary{groupplots}
\usepackage{pgfplotstable}
\usetikzlibrary{spy}
\usepackage[caption=false]{subfig} 
\pgfplotsset{compat=newest}

\newcommand{\setmap}[3]{#1:#2 \mathrel{\vcenter{\mathsurround0pt
\ialign{##\crcr
		\noalign{\nointerlineskip}$\rightarrow$\crcr
		\noalign{\nointerlineskip}$\rightarrow$\crcr
		}}}%
		#3}

\newcommand{\naturals}{\mathbb{N}}
\newcommand{\real}{\mathbb{R}}
\newcommand{\realnonneg}{\mathbb{R}_{\ge 0}}
\newcommand{\realpos}{\mathbb{R}_{> 0}}
\newcommand{\sym}{\operatorname{Sym}}
\newcommand{\until}[1]{{[#1]}}
\newcommand{\map}[3]{#1:#2 \rightarrow #3}
\newcommand{\qedA}{~\hfill \ensuremath{\square}}

\newcommand{\longthmtitle}[1]{\mbox{}{\textit{(#1):}}}
\newcommand{\setdef}[2]{\{#1 \; | \; #2\}}
\newcommand{\setdefb}[2]{\big\{#1 \; | \; #2\big\}}
\newcommand{\setdefB}[2]{\Big\{#1 \; | \; #2\Big\}}
\newcommand*{\SetSuchThat}[1][]{} 
\newcommand*{\MvertSets}{%
    \renewcommand*\SetSuchThat[1][]{%
        \mathclose{}%
        \nonscript\;##1\vert\penalty\relpenalty\nonscript\;%
        \mathopen{}%
    }%
}
\MvertSets 
\DeclarePairedDelimiterX \Set [2] {\lbrace}{\rbrace}
    {\,#1\SetSuchThat[\delimsize]#2\,}

\newtheorem{theorem}{Theorem}[section]

\newtheorem{lemma}[theorem]{Lemma}
\newtheorem{corollary}[theorem]{Corollary}
\newtheorem{remark}[theorem]{Remark}
\newtheorem{definition}[theorem]{Definition}
\newtheorem{assumption}[theorem]{Assumption}

\newtheorem{problem}{Problem}

\renewcommand{\SS}{\mathcal{S}}
\newcommand{\UU}{\mathcal{U}}
\newcommand{\CC}{\mathcal{C}}
\newcommand{\TT}{\mathbf{T}}
\newcommand{\HH}{\mathcal{H}}
\newcommand{\Lie}{\mathcal{L}}
\newcommand{\VV}{\mathcal{V}}
\newcommand{\JJ}{\mathcal{J}}
\newcommand{\xx}{\mathbf{x}}
\newcommand{\uu}{\mathbf{u}}
\newcommand{\vv}{\mathbf{v}}

\newcommand{\kk}{\mathbf{k}}
\newcommand{\ww}{\mathbf{w}}

\newcommand{\BB}{\mathbf{B}}
\newcommand{\cc}{\mathbf{c}}
\newcommand{\WW}{\mathbf{W}}
\newcommand{\PP}{\mathbf{P}}
\newcommand{\sint}{\operatorname{si}}

\newcommand{\ones}{\mathbbm{1}}
\newcommand{\convex}{\text{co}}

\newcommand{\dis}{{\operatorname{dis}}}
\newcommand{\agg}{{\operatorname{agg}}}

\newcommand{\str}{{\operatorname{str}}}
\newcommand{\nom}{{\operatorname{nom}}}
\newcommand{\opt}{{\operatorname{opt}}}
\newcommand{\cm}{{\operatorname{cm}}}
\newcommand{\target}{{\operatorname{target}}}
\newcommand{\Ball}[2]{\mathbb{B}_{#1}(#2)}

\allowdisplaybreaks

\definecolor{matlabBlue}{rgb}{0,0.4470,0.7410}
\definecolor{matlabOrange}{rgb}{0.8500,0.3250,0.0980}
\definecolor{matlabYellow}{rgb}{0.9290,0.6940,0.1250}
\definecolor{matlabPurple}{rgb}{0.4940,0.1840,0.5560}
\definecolor{matlabGreen}{rgb}{0.4660,0.6740,0.1880}
\definecolor{matlabLightBlue}{rgb}{0.3010, 0.7450, 0.9330}
\definecolor{matlabDarkRed}{rgb}{0.6350, 0.0780, 0.1840}
\definecolor{darkGreen}{RGB}{19,122,16}

\pgfplotscreateplotcyclelist{MyCyclelist2}{%
	{matlabBlue, mark = none, line width=0.6pt},
	{matlabOrange, mark = none, line width=0.6pt},
	{matlabYellow, mark = none, line width=0.6pt},
	{matlabPurple, mark = none, line width=0.6pt},
	{matlabGreen, mark = none, line width=0.6pt},
	{matlabLightBlue, mark = none, line width=0.6pt},
	{matlabDarkRed, mark = none, line width=0.6pt},
	{darkGreen, mark = none, line width=0.6pt}
}

\begin{document}
\begin{frontmatter}

  \title{Nonsmooth Control Barrier Function Design of Continuous
    Constraints for Network Connectivity Maintenance} \thanks{This
    work was supported by ONR Award N00014-18-1-2828 and NSF Awards
    ECCS-1917177 and CMMI-2044900. During the preparation of the bulk
    of this work, P. Ong was affiliated with the University of
    California, San Diego.}

\vspace{-10pt}

\author[First]{Pio Ong}%
\author[Second]{\quad Beatrice Capelli}%
\author[Second]{\quad Lorenzo Sabattini}%
\author[Third]{\quad Jorge Cort\'es}

\address[First]{Department of Mechanical and Civil Engineering, California Institute of Technology, USA, {\tt\small
      pioong@caltech.edu}}
\address[Second]{Department of Sciences and Methods for Engineering, University of Modena and Reggio Emilia, Italy,  {\tt\small{\{beatrice.capelli, lorenzo.sabattini\}@unimore.it}}}
\address[Third]{Department of Mechanical and Aerospace Engineering, University of California, San Diego, USA, {\tt\small
      cortes@ucsd.edu}}

\begin{abstract}
This paper considers the problem of maintaining global connectivity of a multi-robot system while executing a desired coordination task.  Our approach builds on optimization-based feedback design formulations, where the nominal cost function and constraints encode desirable control objectives for the resulting input.  Our solution uses the algebraic connectivity of the multi-robot interconnection topology as a control barrier function and critically embraces its nonsmooth nature. We take advantage of the understanding of how Laplacian eigenvalues behave as their multiplicities change, in combination with the flexibility provided by the concept of control barrier function, to carefully design additional constraints that guarantee the resulting optimization-based controller is continuous and maintains network connectivity. The technical treatment combines elements from set-valued theory, nonsmooth analysis, and algebraic graph theory to imbue the proposed constraints with regularity properties so that they can be smoothly combined with other control constraints.  We provide simulations and experimental results illustrating the effectiveness and continuity of the proposed approach in a resource gathering problem.
\end{abstract}
\begin{keyword}
Multi-Robot Systems, Connectivity Maintenance, Nonsmooth Control Barrier Functions, Algebraic Graph Theory
\end{keyword}

\end{frontmatter}

\section{Introduction}

Multi-robot systems can accomplish a variety of tasks through coordinated behavior in many scenarios. Such systems are more versatile, more robust, and better performing than a single specialized robot. To enjoy these advantages, cooperative strategies for multi-robot systems must overcome a number of hurdles, including scalability, graceful degradation with respect to agent failures, and connectivity maintenance, which is the focus of this work. In fact, the ability to interchange information across the network is critical to accomplish emergent coordinated behavior, such as flocking, agreement, coverage, rendezvous, etc., cf.~\citep{FB-JC-SM:09,MM-ME:10,JC-ME:17-jcmsi} and references therein. Connectivity maintenance is hence a fundamental aspect of  cooperative strategies which must be considered in conjunction with the objectives that the multi-robot systems seek to achieve. This integration must be carefully balanced to avoid getting robots stuck in place or display erratic changes in their motions to avoid losing connectivity. Motivated by these observations, this paper investigates how to ensure connectivity while efficiently managing constraints related to the objective of the multi-robot system, with a special emphasis on the continuity of the resulting feedback controller.

\subsubsection*{Literature Review}
Multi-robot systems rely on coordination among agents to achieve their goals. In order to be able to interchange information across the network, the interaction graph must be connected. The concept of algebraic connectivity~\citep{CDG-GFR:01} of a graph, also known as Fiedler eigenvalue~\citep{MF:73}, characterizes the connectivity of a network graph by transforming it into an eigenvalue computation problem. For multi-robot systems, the network graph is dynamically changing as the robots' states evolves and they navigate through their tasks. Typically, robot network graphs are determined via proximity graphs~\citep{FB-JC-SM:09,MZZ-GJP:14-sv}, where the degree of connectivity changes along the robots' trajectories. Connectivity maintenance of dynamic graphs can be categorized into two approaches, local and global, depending on how connectivity is enforced. In the local approach, connectivity is maintained by reasoning over the connections present in the initial graph. This includes the direct method of preserving all initial connections, \citep[see e.g.,][]{MJ-ME:07}, which limits the graph to one arrangement. This method can be improved by considering instead multiple-hops neighbors and allowing rearrangements in the edges~\citep{MMZ-GJP:05,MDS-JC:09-sicon}, but its flexibility is still limited by the initial robot configuration. The global approach reasons more broadly over network connectivity using network-wide metrics such as algebraic connectivity.  Under this approach, we find works that pose connectivity as a problem of maximizing algebraic connectivity~\citep{SB:06,YK-MM:06}. The idea is to find a robot motion that will increase the algebraic connectivity. A decentralized implementation of this idea is explored in \citep{MDG-AJ:06}. Nevertheless, maximizing the algebraic connectivity in all scenarios can be overly restrictive. In this regard,~\citep{LS-NC-CS:13,MDS-JC:09} introduce more flexibility by allowing algebraic connectivity to decrease when its value is large.

Our connectivity maintenance solution here is based on the concept of Control Barrier Function (CBF) from the safety-critical control literature. Control Barrier Functions~\citep{PW-FA:07} build on the barrier certificate~\citep{SP-AJ:04} notion, and is used to find choices of control inputs that makes the certificate increase, guaranteeing forward invariance of a desired set. The CBF idea can be refined further by abandoning the monotonicity of the certificate. This idea is related to the concept of practical stability with Lyapunov functions, with an additional restriction on the evolution of the certificate within the desired set~\citep{PO-AB-TH-LK-PS:06}. It is later formalized in the context of safety~\citep{ADA-SC-ME-GN-KS-PT:19} by using Nagumo theorem~\citep{FB-SM:07} as the basis for set invariance. This refined version introduces the concept of letting the certificate also decrease depending on the level of safety. In the context of connectivity maintenance, CBFs flexibly allows algebraic connectivity to decrease as long as the graph does not become disconnected. CBFs are employed in both aforementioned connectivity maintenance approaches in \citep{ME-JNP-GN-SH:18} and \citep{BC-LS:20}, respectively. Regarding the latter, there is no guarantee on the continuity of the proposed feedback controller because of the lack of smoothness of the algebraic connectivity. Here instead, we rely on Nonsmooth Control Barrier Functions (NCBF)~\citep{PG-JC-ME:17}, a generalization of CBF, to properly account for the nonsmoothness of algebraic connectivity and ensure the continuity of the resulting feedback controller.

Controllers that utilize CBFs are typically based on optimization formulations, \citep[see e.g.,][]{ADA-XX-JWG-PT:17,ADA-SC-ME-GN-KS-PT:19}. For this type of controllers, there are multiple approaches to determine continuity. Using perturbation theory, the paper~\citep{BJM-MJP-ADA:15} studies the smoothness properties of optimization-based controllers with CBFs but the result is only applicable to continuously differentiable CBFs. From a set-valued theory perspective, \citep{RAF-PVK:96} shows continuity of minimum-norm controllers, i.e., when the objective function is a norm. For more general objective functions (like the one considered here), we rely on Berge Maximum Theorem~\citep{CDA-KCB:99}, a well-known result in parametric optimization, to guarantee continuity of the feedback controller. 

\subsubsection*{Statement of Contributions}
This paper considers a multi-robot system with fully actuated first-order dynamics. The underlying interaction network is described by a continuously differentiable proximity graph. We address the problem of maintaining global connectivity of the multi-robot system that is operating under some nominal control constraints. The contributions of the paper are threefold. The first contribution is the synthesis of two different set-valued constraint maps for global connectivity maintenance. The proposed constraints are based on  NCBFs and are able to handle, in a continuous way,  the abrupt changes caused by the jumps in multiplicity of the algebraic connectivity as a function of the network state. Establishing this fact relies on a careful application of various notions and results from set-valued analysis. With our proposed constraints, the resulting optimization-based controller is continuous, which thereby guarantees the existence of a solution and avoids issues such as chattering in its discrete-time implementation.
Our second contribution deals with the well-posedness of the considered problem. As we allow for the possibility of the network to have control constraints beyond connectivity maintenance, one question that we answer is in regard to the existence of a solution to our problem, i.e., a continuous controller that can both respect control constraints and maintain network connectivity. We use a generalization of Artstein's theorem \citep{ZA:83}
to deduce a mild and verifiable condition that guarantees our problem is well-posed. Our final  contribution are the continuity results, as a  function of the network state, for any intersection of eigenspaces of the graph Laplacian. We rely on this result to study how the algebraic connectivity changes. Since the generalized gradient of algebraic connectivity is related to its associated eigenspace, the results add to the literature on regularity of algebraic connectivity. We believe our second and third contributions may have useful applications beyond the subject matter of this paper. We conclude the paper by illustrating the effectiveness of our results in a resource gathering problem, both in simulations and an experiment. The problem consists in a group of robots trying to reach their assigned target locations which cannot do so without losing connectivity. We show that using our proposed results, each robot in the network completes its tasks with continuous feedback control inputs, and the robot network remains connected throughout.

A preliminary version of this paper appeared at the IEEE Conference on Decision and Control~\citep{PO-BC-LS-JC:21}. The added value of the present work is justified by the following additions: (i) a more general control synthesis problem formulation that incorporates a nominal constraint map, which results in a more challenging technical analysis; (ii) the generalization of Artstein's theorem to formulate a reasonable assumption for the well-posedness of the problem; (iii) the establishment of the continuity property of merged eigenspaces of the graph Laplacian as a function of the network state, which was only speculated in the preliminary version of the paper; (iv) the new simulation example along with a validation of the results in an experiment on four small-wheeled robots. In addition, we provide throughout the paper all the necessary background and discussions on intuitions behind the proposed ideas.

\section{Preliminaries}
This section introduces basic notation and key concepts from graph theory, set-valued and nonsmooth analysis, and Nonsmooth Control Barrier Functions.

\subsection{Notation} 
The symbols $\naturals$, $\real$, $\realnonneg$, and $\realpos$ represent the set of natural, real, real nonnegative, and real positive numbers, respectively. We write $\sym^n$ for the space of $n\times n$ symmetric matrices with real values. For $m,n \in \naturals$, we denote $\until{m:n} = \{m,\dots,n\}$, and we write $\until{1:n}$ simply as $\until{n}$. Given a finite set $\mathcal I$, $\vert \mathcal I \vert$
is its cardinality. The convex closure of a set~$\SS$ is represented by $\convex (\SS)$. Given $\xx\in\real^N$, $\|\xx\|$ denotes its Euclidean norm. We use the symbol $\ones$ for the vector of all ones (of appropriate dimension). The unit sphere in $\real^n$ is denoted by $\mathbb S^n = \setdefb{\vv\in \real^n}{\|\vv\| = 1}$. The open ball of radius $\delta>0$ centered at $\xx^*\in\real^N$ is $\Ball{\delta}{\xx^*} = \setdefb{\xx\in\real^N}{\|\xx-\xx^*\|<\delta}$. Given matrices $\mathbf{A},\BB\in \real^{n\times n}$, the Frobenius product is $\mathbf{A}\cdot \BB = \sum_{i,j} \mathbf{A}_{ij}\BB_{ij}$. We note the property that $\vv \vv^\top\cdot \mathbf{A} = \vv^\top \mathbf{A}\vv$, for $\vv\in \real^n$. The Frobenius norm is given by $\|\mathbf{A}\|_F=(\mathbf{A}\cdot \mathbf{A})^{1/2}$. A continuous function $\map{\alpha}{\real}{\real}$ is of extended class $\mathcal{K}$ if $\alpha$ is strictly increasing, and $\alpha(0) = 0$.
Moreover, $\text{supp(f)}$ is the support of the function~$f$, i.e., the set of $\xx$ where $f(\xx)\neq 0$.

\subsection{Graphs and Laplacian Spectrum}

A graph is a triplet $\mathcal{G}= (V, E, \mathbf A)$, where $V$ is a set of vertices,  $E \subseteq  V \times  V$ is a set of edges, and $\mathbf A\in \real^{\vert  V\vert \times \vert V\vert}$ is the adjacency matrix, with $\mathbf A_{ij}>0$ if $(i,j)\in E$, and  $\mathbf A_{ij}=0$ otherwise. 
The graph is undirected if $\mathbf A$ is symmetric. A path is an ordered sequence of vertices such that all pairs of consecutive vertices are elements of~$E$. The graph is connected if there exists a path between any two vertices. The degree matrix $\mathbf D\in \real^{\vert  V\vert \times \vert V\vert}$ is a diagonal matrix whose $i$th element is $\mathbf \mathbf{D}_{ii} = \sum_{j \in  V}\mathbf A_{ij}$.
The Laplacian matrix $\mathbf L$, defined by $\mathbf L := \mathbf D - \mathbf A$, is symmetric and positive semidefinite, and consequently has real and nonnegative eigenvalues. We denote these eigenvalues with $\phi_m \in \realnonneg$, ordering them in an increasing manner with the subscripts $m\in \until{\vert V\vert}$, i.e., $0=\phi_1 \leq \phi_2 \leq \dotsc \leq \phi_{\vert V\vert}$. The eigenvalue $\phi_1=0$ is simple (with associated eigenvector $\ones$) if and only if the graph is connected. This justifies the terminology of $\phi_2$ as the algebraic connectivity (also known as Fiedler eigenvalue). For network systems, graphs are used to described the underlying interaction topology, and they can vary according to the system states. A state-dependent graph $\xx \mapsto \mathcal{G}(\xx)$ is called a \emph{proximity graph}~\citep{FB-JC-SM:09}. In such a case, the Laplacian matrix $\xx \mapsto \mathbf{L}(\xx)$ is then also a function of the state. We define the function $\lambda_m(\xx):= (\phi_m \circ \mathbf L)(\xx)$ to be the Laplacian's eigenvalues as a function of the state. Given a trajectory $t\mapsto \xx(t)$, a graph remains robustly connected at all times if $\lambda_2(\xx(t))\geq \varepsilon$, where $\varepsilon \in \realpos$ is a threshold parameter providing a robustness margin in ensuring connectivity. 

\subsection{Continuity of Set-Valued Maps}
A set-valued map $\setmap{\UU}{\real^N}{\real^M}$ assigns a subset of $\real^M$ to each point in $\real^N$. A set-valued map $\UU$ is closed-valued, convex-valued, compact-valued, and has a nonempty interior if its image at each point of its domain is closed, convex, compact, and has a nonempty interior, respectively. All set operations,  e.g., union and intersection, between set-valued maps are performed pointwise.
Throughout the paper, we consider set-valued maps arising from a single-valued function $\map{g}{\real^N\times\real^M}{\real^d}$ as follows:
\begin{equation}
\label{eq:generic}
\UU(\xx) = \setdef{\uu\in\real^M}{g(\xx,\uu)\leq 0}.
\end{equation}
Given $\xx$, we say $\uu$ strictly satisfies $\UU(\xx)$ if $g(\xx,\uu)<0$.

The concept of continuity for set-valued maps is more intricate than the one for single-valued functions. Continuity of set-valued maps is often broken down into different types of hemicontinuity. Here we present the two that we rely on: upper and lower hemicontinuity\footnote{Sometimes referred to as semicontinuity, see e.g., \citep{AL-AS:85}.}.

\begin{definition}\longthmtitle{Set-Valued Map Continuity~\citep{KCB:85}}
A set-valued map $\setmap{\UU}{\real^N}{\real^M}$ is
\begin{itemize}
    \item upper hemicontinuous (UHC) at $\xx$ if for any neighborhood $\bar \UU$ of $\UU(\xx)$, there exists $\delta>0$ such that, if $\|\xx-\xx'\|<\delta$, then $\UU(\xx') \subset \bar \UU$;
    \item lower hemicontinuous (LHC) at $\xx$ if for each $\uu \in \UU(\xx)$ and for any sequence $\{\mathbf{x}^k\}_{k\in \naturals}$ converging to $\xx$, there exists a sequence $\{\uu^k\}_{k\in \naturals}$ converging to $\uu$ with $\uu^k \in \UU(\xx^k)$;
    \item continuous at $\xx$ if it is both UHC and LHC at $\xx$.
\end{itemize}
\end{definition}
Note here that UHC and LHC are equivalent for single-valued functions. For convenience, the map is (hemi)continuous if it is (hemi)continuous for all $\xx$. Interestingly for set-valued maps of the form~\eqref{eq:generic}, even $g$ being continuous is not enough to ensure the map $\UU$ is continuous. In fact, to ensure UHC and LHC, we will resort to the additional requirements stated in the following results.

\begin{lemma}
\label{lem:UHC_req}
\longthmtitle{UHC Requirements \citep[Lem 5.7]{GS:18}}
Assume $g$ is continuous.
If $g$ is convex in $\uu$, and $\UU(\xx)$ is nonempty and compact at $\xx$, then $\UU$ is UHC at $\xx$.~\hfill $\qedA$
\end{lemma}

\begin{lemma}
\label{lem:LHC_req}
\longthmtitle{LHC Requirements \citep[Lem 5.2]{GS:18}}
Assume $g$ is  continuous. If $\UU$ has a nonempty interior and is convex-valued, then $\UU$ is LHC.~\hfill $\qedA$
\end{lemma}

In our treatment, we also rely on various results on how hemicontinuity is preserved under set-valued map intersections.

\begin{lemma}\longthmtitle{Intersection of UHC maps \citep[11.21a]{KCB:85}}
\label{lem:insct-uhc}
Let the set-valued maps $\setmap{\UU_1,\UU_2}{\real^N}{\real^M}$ be UHC and closed-valued at $\xx$. The intersection $\UU_1 \cap \UU_2$ is also UHC at $\xx$ if it is nonempty at $\xx$.~\hfill \qedA
\end{lemma}

\begin{lemma}\longthmtitle{Intersection of LHC maps \citep[Thm. B]{AL-AS:85})}
\label{lem:insct-lhc}
Let the set-valued maps $\setmap{\UU_1,\UU_2}{\real^N}{\real^M}$ be LHC and locally convex-valued at $\xx$. The intersection $\UU_1 \cap \UU_2$ is also LHC at $\xx$ if it has a nonempty interior at $\xx$.~\hfill \qedA
\end{lemma}

\subsection{Nonsmooth Analysis}
Here we present basic notions of nonsmooth analysis following~\citep{FHC:83}. 
Given a locally Lipschitz function $\map{h}{\real^N}{\real}$, the generalized directional derivative of $h$ at $\xx \in \real^N$ in the direction $\mathbf d \in \real^N$ is
$$
h^\circ(\xx;\mathbf d) = \limsup_{\xx'\rightarrow \xx,s\downarrow 0}\frac{h(\xx'+s\mathbf d)-h(\xx')}{s}.
$$
The generalized gradient of $h$ at $\xx$ is then given by
$$
\partial h(\xx) = \setdef{\zeta\in \real^N}{h^\circ(\xx;\mathbf d)\geq \zeta^\top\mathbf d,~\forall \mathbf d\in \real^N}.
$$
If the function $h$ is continuously differentiable at $\xx$, the generalized gradient is a singleton, $\partial h(\xx) = \{\nabla h(\xx)\}$.

In our analysis, we find it useful to describe how a nonsmooth function changes along the trajectories of a dynamical system. Consider the nonlinear system,
\begin{equation}\label{sys:gen}
    \dot \xx = f(\xx,\uu) ,
\end{equation}
with $\map{f}{\real^N\times \real^M}{\real^N}$ measurable and essentially locally bounded,
where $\xx$ is the state and $\uu$ is the control input.
The weak set-valued Lie derivative~\citep{PG-JC-ME:17,DS-BP:94} is
$$
\Lie_fh(\xx,\uu) = \setdefb{\zeta^\top f\in\real}{\zeta\in\partial h(\xx)}.
$$
The Lie derivative describes the rate of change of $h$ along a trajectory of the system. Let $t\rightarrow\uu (t)$ be a control signal, and $t\rightarrow\xx(t)$ be a Carath\'eodory solution\footnote{A Carath\'eodory solution is an absolutely continuous trajectory that satisfies the system dynamics at almost every time, in the sense of Lebesgue measure.} to the differential equation~\eqref{sys:gen}, then 
\begin{align}\label{eq:rate-of-change}
\frac{d}{dt}h(\xx(t))\in \Lie_fh(\xx(t),\uu(t)),~a.e.~t\geq 0.
\end{align}
In essence, the weak set-valued Lie derivative contains all the possible rates of change of the function $h$ along a solution of the dynamical system.

\subsection{Nonsmooth Control Barrier Functions}\label{sec:cbf}

We use Nonsmooth Control Barrier Functions (NCBF)~\citep{PG-JC-ME:17} to establish forward invariance of a desired set. Consider the dynamical system~\eqref{sys:gen} and a set
$\CC = \setdef{\xx\in\real^N}{h(\xx)\geq 0}$ with a locally Lipschitz continuous $\map{h}{\real^N}{\real}$, referred to as a nonsmooth control barrier function. Indeed, for a continuous trajectory $t\rightarrow \xx(t)$, we can ensure $h$ remains positive if we constrain $h$ from decreasing whenever $h(\xx(t))=0$. This can be done by imposing a constraint, as a function of network state $\xx$, on our choice of the input $\uu$ with a set-valued map
$$
\UU(\xx) = \setdefB{\uu \in \real^M}{\min\Lie_Fh(\xx,\uu)\geq -\alpha(h(\xx))} ,
$$
where $\alpha$ is a locally Lipschitz extended class $\mathcal{K}$ function. Given~\eqref{eq:rate-of-change}, by taking the minimum element of the set-valued Lie derivative, the constraint map enforces the bound even for the worst-case rate of change of $h$. Note importantly that the above constraint map does not only limit the choice of $\uu$ for $\xx$ at the boundary of $\CC$ where $h(\xx)=0$, but also in the interior where $h(\xx)>0$, even when it is not necessary. Rather than outright allowing any choice of $\uu$, the constraint map gradually becomes stricter for states closer to the boundary. The idea here is to begin consider the necessary constraint as the trajectory approaches the boundary, and thereby provide some robustness to how the set $\CC$ is rendered forward invariant.

\section{Problem Statement}\label{sec:problem_statement}
Consider a group of $n$ robots, evolving according to a single-integrator dynamics of the form
\begin{equation}
\label{sys:full}
	\dot{x}_r = u_r, \quad \forall r\in\until{n} ,
\end{equation}
where $x_r\in\real^{d_r}$ and $u_r \in \real^{d_r}$ are the state and the control input associated with the $r$-th robot (note that the state dimensions of each robot might be different). For convenience, we define state and input variables for the network system as follows: let $N=\sum_{r\in\until n} d_r$ and denote $\xx=\left[ x_1^\top, \dotsc, x_n^\top \right]^\top \in \real^{N}$ and $\uu=\left[ u_1^\top, \dotsc, u_n^\top \right]^\top \in \real^{N}$. We use the shorthand notation $\map{f_{\sint}}{\real^N \times \real^N}{\real^N}$ to refer compactly to the dynamics~\eqref{sys:full} for the whole group of agents. The underlying interaction topology is specified by a proximity graph $\xx \mapsto \mathcal{G}(\xx)=(\until n, E(\xx),\mathbf{A}(\xx))$, for which we assume that the function $\xx \mapsto \mathbf{A}(x)$ is continuously differentiable\footnote{This assumption is satisfied by commonly employed weight assignments~\citep{MDS-JC:09, AG-LS-GU:17}.}.

We are interested in designing a continuous controller $\map{\kk}{\real^N}{\real^N}$ such that the network system under feedback $\uu=\kk(\xx)$ enjoys some desirable performances and asymptotic guarantees. Continuity is an important property, both from a theoretical and practical viewpoint. 
Regarding the former, continuity guarantees the existence of Carath\'eodory (in fact, classical) solutions~\citep[Thm. 5.1]{JKH:69}. At the same time, continuity makes it easier for the desired feedback control signal to be implemented on digital platforms.

A commonly used design methodology to synthesize controllers is based on optimization and takes the form
\begin{align}\label{eq:controller}
\kk_\opt(\xx) = \underset{\uu\in \UU(x)}{\arg\!\min}~J(\xx,\uu) ,
\end{align}
where $\map{J}{\real^N \times \real^N}{\real}$
is a cost function encoding some desirable objective (e.g., minimal deviation from a prescribed input, minimum-energy control specifications) and $\setmap{\UU}{\real^N}{\real^N}$ is a set-valued map encoding constraints on the control input at each~$\xx$ (e.g., bounds on magnitude, stability performance using control Lyapunov function). This formulation is flexible as it can address simultaneously different performance requirements: the map $\UU$ can be itself an intersection of multiple set-valued maps, each representing a different control constraint from a performance aspect (input boundedness, infinitesimal decrease of certificate).

We consider the scenario where the robot group has a nominal control constraint map $\xx \mapsto \UU_{\nom}(\xx)$, defined via a function $\map{g_\nom}{\real^N\times\real^N}{\real^{d_\nom}}$ as
$$
\UU_\nom(\xx) = \setdef{\uu \in \real^N}{g_\nom(\xx,\uu)\leq 0}.
$$
The components of $g_\nom$ here represent constraints that the robot group must respect to achieve different control performances and goals. This nominal constraint map, however, does not encode any network connectivity constraint. We are then interested in solving the following problem.

\begin{problem}\longthmtitle{Continuous Connectivity Controller Design Problem}
\label{prob:connect}
{\rm
Consider the multi-robot system~\eqref{sys:full} operating with the optimization-based controller~\eqref{eq:controller}. Design the constraint map~$\UU$ so that:
\begin{itemize}
    \item the controller $\kk_\opt$ is continuous;
    \item the nominal constraint map is respected, i.e., $\UU\subseteq\UU_\nom $;
    \item the underlying graph $\mathcal G$ remains connected at all time.~\hfill $\bullet$
\end{itemize}}
\end{problem}

We make the following assumptions on the cost function $J$ and the nominal constraint map~$\UU_\nom$ to make sure Problem~\ref{prob:connect} is solvable. First,~$\UU_\nom$ should be large enough so that, at  each  state, there exists a control  that  can simultaneously maintain connectivity and satisfy  the  nominal  constraints (we formulate this assumption mathematically later in our technical discussion, cf. Remark~\ref{rmk:feasible}). As one can expect, continuity of $\kk_\opt$ is related to continuity of the cost function $J$ and the constraint map~$\UU$. In this regard, Berge Maximum Theorem~\citep[Thm. 17.31]{CDA-KCB:99} states that, if $J$ and $\UU$ are continuous, $\UU$ is compact-valued, and the resulting $\kk_\opt$ is single-valued, then $\kk_\opt$ is continuous. Based on this result, we make the following continuity assumption.

\begin{assumption}
\label{assump:continuity}
\longthmtitle{Continuity Assumption on Cost and Nominal Constraint}
The functions $J$ and $g_\nom$ are continuous.~\hfill $\bullet$ 
\end{assumption}

We do not make a direct assumption on the continuity of~$\UU_\nom$ for greater generality. In fact, such assumption would rule out many commonly used constraint maps  (e.g., control affine constraint maps are typically not UHC). As such, we rely instead on the following assumption.

\begin{assumption}\longthmtitle{Convexity Assumption on Cost and Nominal Constraint}\label{assump:convexity}
The function $J$ is strictly convex in $\uu$ and $g_\nom$ is convex in $\uu$.~\hfill $\bullet$
\end{assumption}

Although convexity is not required by Berge Maximum Theorem, the above assumption is justified by several reasons. First, the assumption helps us make the optimization problem that defines the controller a convex program, which opens the way to employing available convex optimization methods to compute the controller. In addition, the strict convexity assumption also ensures that the controller is single-valued for any given $\xx$, which is a requirement of Berge Maximum Theorem. More importantly, the convexity assumption also opens up the possibility of $\UU$ being defined by unbounded constraints, despite the compact-valued requirement in Berge Maximum Theorem. To reconcile this, we  consider the sublevel sets of~$J$. Suppose for each $\xx$, there exists a control $\xx\mapsto \bar \uu(\xx)$ such that $\bar \uu(\xx) \in \UU(\xx)$, and define
\begin{equation}\label{eq:J_sublvl}
\JJ_{\bar \uu}(\xx) = \setdef{\uu\in \real^{N}}{\|J(\xx,\uu)\|\leq \|J(\xx,\bar \uu(\xx))\|+\delta_J}
\end{equation}
with $\delta_J \in \realpos$. Note that this set-valued map is compact-valued due to strict convexity of~$J$.  In addition, when $\JJ_{\bar \uu}$ is considered in conjunction with $\UU$, it is always inactive at the optimizer because $\bar \uu$ is a feasible point. Consequently, for a properly designed $\UU$, even if it is not compact-valued, we may consider $\UU\cap \JJ_{\bar \uu}$ as the constraint map without changing the optimizer at each $\xx$ and apply Berge Maximum Theorem.

\section{Discontinuity in the Naive Connectivity Maintenance Solution}\label{sec:naive_solution}

In this section we make a first attempt at 
solving Problem~\ref{prob:connect} using algebraic connectivity as a nonsmooth control barrier function. We show that the proposed solution falls short because the resulting feedback controller is discontinuous. This exercise serves two purposes. On the one hand, it motivates the technical refinement pursued in our exposition later. On the other, it helps us pinpoint the obstructions associated with solving Problem~\ref{prob:connect}, providing the necessary exposition for the rationale behind our solutions.

For maintaining connectivity, it seems natural to use the algebraic connectivity as a NCBF to guarantee $\lambda_2$ remains positive along the trajectory. This is essentially the approach taken in~\citep{BC-LS:20} (with the difference that we explicitly account for the nonsmoothness of $\lambda_2$ in the exposition here).
Consider the safe set of connected robot configurations
\[
\CC_\varepsilon := \setdefb{\xx\in\real^N}{\lambda_2(\xx) \geq \varepsilon} ,
\]
with $\varepsilon\in\realpos$. The parameter $\varepsilon$ is introduced here to provide robustness and ensure the safe set is closed. Its introduction makes the safe set smaller and hence the safety specification more conservative.  Let $\xx \mapsto h(\xx) = \lambda_2(\xx) - \varepsilon$ be our candidate NCBF. Resorting to the discussion of Section~\ref{sec:cbf}, we specify a constraint map for the purpose of connectivity maintenance.

\begin{lemma}\longthmtitle{Connectivity Maintenance Constraint Map}\label{lem:connect}
Consider the multi-robot system~\eqref{sys:full} operating with a controller $\xx \mapsto \kk(\xx)$. Given a locally Lipschitz extended class $\mathcal K$ function~$\alpha$, define the constraint map
\[
\UU_\cm(\xx) := \setdefb{\uu\in \real^N}{
\min \Lie_{f_{\sint}}\lambda_2(\xx,\uu) \geq -\alpha(\lambda_2(\xx)-\varepsilon)}.
\]
If $\kk(\xx)\in \UU_\cm(\xx)$ for all $\xx\in\CC_\varepsilon$, then for any initial connected network configuration $\xx_0 \in \CC_\varepsilon$, $\lambda_2(\xx(t)) \ge \varepsilon$ along all Carath\'eodory solutions of the closed-loop system under 
$\uu=\kk(\xx)$, ensuring that network connectivity is maintained.
\end{lemma}

Lemma~\ref{lem:connect} is a direct result of using $h(\xx) = \lambda_2(\xx)-\varepsilon$ as a NCBF, cf. \citep[Thm.3]{PG-JC-ME:17}. Our first attempt to utilize the result is to design an optimization-based controller~\eqref{eq:controller} naively defined with the connectivity maintenance constraint map,
\begin{align}\label{eq:dis_controller}
\kk_{\dis}(\xx) := \underset{u\in \UU_\cm(\xx)}{\arg\!\min}~J(\xx,\uu).
\end{align}
Unfortunately, this controller is not guaranteed to be continuous. This can cause a number of undesired phenomena. For example, sample-and-hold implementations of the controller may exhibit chattering behavior because its continuous-time counterpart is not continuous (cf., Fig. \ref{fig:simulation_input_3} in Sec. \ref{sec:experiment}). More importantly for our problem, the connectivity maintenance result provided by Lemma~\ref{lem:connect} is only guaranteed along Carath\'eodory solutions, and such solutions might not exist if the controller is not continuous. As such, discontinuous controllers like $\kk_{\dis}$ may fail to maintain connectivity of the multi-robot system.

The discontinuity issue arises because $\UU_\cm$ itself is not continuous and does not meet the requirement of Berge Maximum Theorem. 
To pinpoint the root cause of the discontinuity of $\UU_\cm$, we review the generalized gradient of the Laplacian's eigenvalues.
Each eigenvalue function $\phi_m$ is globally Lipschitz with respect to the entries of the Laplacian matrix (cf.,~\citep[Lem. 1]{MDS-JC:09} and~\citep[Thm. 2.4]{ASL:96}). As a result, if  $\mathbf{L}$ is a continuously differentiable function of the network state, then $\lambda_m = \phi_m \circ \mathbf{L}$ is also Lipschitz. Therefore, generalized gradients are well-defined for the eigenvalue functions. Mathematically, the generalized gradient of $\phi$ is given by, cf.~\citep[Thm. 1]{MDS-JC:09},
\begin{equation}
\label{eq:gen_grad}
 \partial \phi_m(\mathbf{L}) = \convex\setdefb{\vv_m\vv_m^\top}{\vv_m \in \VV_m(\mathbf{L})},
\end{equation}
where $\VV_m(\mathbf{L}) := \setdef{\vv_m\in \mathbb S^n}{\mathbf{L}\vv_m = \phi_m(\mathbf{L}) \vv_m}$ is the set of normalized eigenvectors associated with $\phi_m$. Using the nonsmooth chain rule~\citep[Thm. 2.3.10]{FHC:83}, the expression for the weak set-valued Lie derivative~\citep[Rmk. 2.1]{PG-JC-ME:17} of $\lambda_m$ with respect to the system~\eqref{sys:full} is 
$$
\Lie_{f_{\sint}}\lambda_m(\xx,\uu) = \partial\phi_m(\mathbf{L}(\xx)) \cdot \Big(\sum_{i\in\until{N}} \frac{\partial \mathbf{L}}{\partial \xx_i} \uu_i \Big).
$$
In the constraint map $\UU_\cm$, we use the minimal value of this set to bound the worst-case rate of change of $\lambda_m$ along the control choice $\uu$. Unfortunately, this minimal value is not a continuous function of~$\xx$. The following result helps us understand why.

\begin{lemma}\longthmtitle{Equivalent Minimization of the Eigenvalue's Set-Valued Lie Derivative}\label{lem:equiv_min} Consider the multi-robot system \eqref{sys:full}. For $m \in \until{N}$, 
let $(\xx,\uu) \mapsto \mu_m(\xx,\uu)$, 
\begin{equation}\label{eq:equiv_min}
    \mu_m(\xx,\uu) := \min_{\vv\in\VV_{m}(\mathbf{L}(\xx))}\vv^\top \Big(\sum_{i\in\until{N}} \frac{\partial \mathbf L}{\partial \xx_i} \uu_i \Big)\vv.
\end{equation}
Then $\min \Lie_{f_{\sint}}\lambda_m(\xx,\uu) = \mu_m(\xx,\uu)$ for any $\xx$ and~$\uu$.
\end{lemma}
\begin{pf}
Let $\PP\in \partial\phi_m(\mathbf{L}(\xx))$ be the element of the generalized gradient~\eqref{eq:gen_grad} corresponding to the minimum value in $\Lie_{f_{\sint}}\lambda_m(\xx,\uu)$, i.e.,
$$
\min \Lie_{f_{\sint}}\lambda_m(\xx,\uu) = \PP\cdot \Big(\sum_{i\in\until{N}} \frac{\partial \mathbf L}{\partial \xx_i} \uu_i \Big).
$$
From the Carath\'eodory theorem of convex hulls~\citep[Thm. 17.1]{RTR:70}\footnote{ The theorem is applicable to the matrix space $\real^{n\times n}$ because its vectorized version is the $n^2$-dimensional real space $\real^{n^2}$.}, since $\PP\in \real^{n\times n}$, there exists $n^2+1$ (not necessarily distinct) points $\{\PP_j\}_{j=1}^{n^2+1}$, each belonging to the set $\partial \phi_m'(\mathbf{L}(\xx)) = \setdefb{\vv_m\vv_m^\top}{\vv_m\in\VV_m(\mathbf L(\xx))}$, such that $\PP$ is a convex combination of $\{\PP_j\}_{j=1}^{n^2+1}$. That is, there exist $0\leq\sigma_j\leq 1$ with $\sum_{j=1}^{n^2+1}\sigma_j=1$ such that
\begin{align}\label{eq:aux}
\min \Lie_{f_{\sint}}\lambda_m(\xx,\uu) &= (\sum_{j=1}^{n^2+1} \sigma_j\PP_j)\cdot \Big(\sum_{i\in\until{N}} \frac{\partial \mathbf L}{\partial \xx_i} \uu_i \Big)\\
&=\sum_{j=1}^{n^2+1} \sigma_j\underbrace{\bigg(\PP_j\cdot \Big(\sum_{i\in\until{N}} \frac{\partial \mathbf L}{\partial \xx_i} \uu_i \Big)\bigg)}_{\geq \min \Lie_{f_{\sint}}\lambda_m(\xx,\uu)}. \notag
\end{align}
From this, there must exist at least one $\PP_j$ such that 
$$
\min \Lie_{f_{\sint}}\lambda_m(\xx,\uu) = \PP_j\cdot \Big(\sum_{i\in\until{N}} \frac{\partial \mathbf L}{\partial \xx_i} \uu_i \Big).
$$
(Otherwise, if $\PP_j\cdot \Big(\sum_{i\in\until{N}} \frac{\partial \mathbf L}{\partial \xx_i} \uu_i \Big)>\min \Lie_{f_{\sint}}\lambda_m(\xx,\uu)$ for every $j \in \until{n^2+1}$,
one would reach from~\eqref{eq:aux}  and $\sum_{j=1}^{n^2+1}\sigma_j=1$ the contradiction that $\min \Lie_{f_{\sint}}\lambda_m(\xx,\uu)>\min \Lie_{f_{\sint}}\lambda_m(\xx,\uu)$).
Hence, it is equivalent to compute the minimization of the Lie derivative on the set $\partial \phi_m'$ rather than its convex closure~$\partial \phi_m$, that is,
\begin{align*}
    \min \Lie_{f_{\sint}}\lambda_m(\xx,\uu) &=
    \min
    \partial\phi_m'(\mathbf{L}(\xx)) \cdot \Big(\sum_{i\in\until{N}} \frac{\partial \mathbf{L}}{\partial \xx_i} \uu_i \Big)
    \\
    &= \min_{\PP'\in\partial \phi_m'(\mathbf L(\xx))} \PP' \cdot \Big(\sum_{i\in\until{N}} \frac{\partial \mathbf L}{\partial \xx_i} \uu_i \Big)\\
    &= \min_{\vv\in\VV_{m}(\mathbf{L}(\xx))} \vv\vv^\top \cdot \Big(\sum_{i\in\until{N}} \frac{\partial \mathbf L}{\partial \xx_i} \uu_i \Big)\\
    &= \mu_m(\xx,\uu),
\end{align*}
where we have used a property of the Frobenius product in the last step. This concludes the proof.~\hfill~\qed
\end{pf}

\begin{figure}[tb]
    \centering
    \includegraphics[width=\columnwidth]{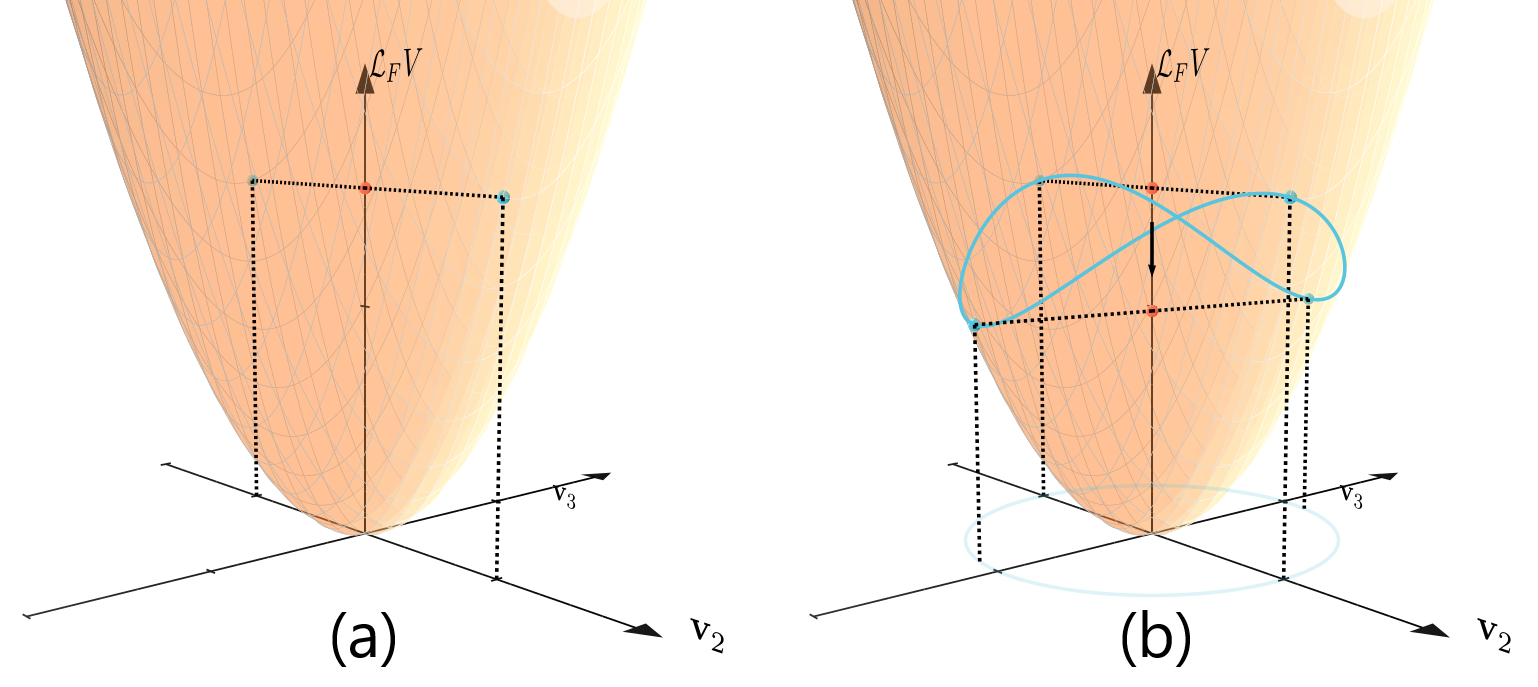}
    \caption{A jump in the minimum value of the set-valued Lie derivative may occur when the corresponding eigenspace expands. Plot (a) shows the optimal value of a quadratic function over a  normalized eigenspace in one dimension (which consists simply of two points, $\vv_2$ and $-\vv_2$). Plot (b) shows that the optimal value of the quadratic function abruptly drops when the normalized eigenspace increases by one dimension to a circle.}
    \label{fig:multiplicity}
\end{figure}

Lemma~\ref{lem:equiv_min} transforms the minimization of the set-valued Lie derivative into an equivalent one with respect to eigenvectors. From this perspective, it is easy to identify the reason for the discontinuity in the minimum value. Whenever the multiplicity of an eigenvalue changes, so does the  dimension of its eigenspace. Consequently, the minimization may abruptly change in value simply because of the abrupt change in the minimization constraint, as illustrated in Fig.~\ref{fig:multiplicity}. We rely on this key insight to synthesize our design in the next section.

\section{Continuous Connectivity Maintenance Constraint Maps}\label{sec:results}
In this section, we propose our solution to Problem~\ref{prob:connect}. We construct two constraint maps for the purpose of connectivity maintenance. The first solution directly addresses the discontinuity issue in the naive solution. This is done by adjusting conservatively the discontinuous term discussed in Section~\ref{sec:naive_solution}. Our second solution refines the first to reduce its conservatism. For clarity of exposition, here we just explain the proposed solutions, and delay the formal technical analysis to Section~\ref{sec:analyses} below.

We first design a connectivity maintenance constraint map by replacing the discontinuous term~$\mu_m$. The discontinuity in $\mu_m$ is due to the abrupt change in the eigenspace being considered in the minimization~\eqref{eq:equiv_min}. One possible fix is to augment the eigenspace preemptively so that there is no abrupt expansion. For $\mathcal I \subseteq \until{n}$, consider
$$
\VV_{\mathcal I} (\xx) := \text{span}\Big\{\bigcup_{p \in \mathcal I}\VV_p(\xx)\Big\}\cap \mathbb S^n,
$$
the normalized span of eigenspaces corresponding to the eigenvalues $\{ \lambda_p \}_{p\in \mathcal I}$ at~$\xx$. We refer to the set-valued map $\VV_{\mathcal I}$ as the \textit{normalized merged eigenspace}. We use this set-valued map to define
\begin{align}\label{eq:mu-I}
\mu_{\mathcal I}(\xx,\uu) := \min_{\vv\in\VV_{\mathcal I}(\xx)}\vv^\top \Big(\sum_{i\in\until{N}} \frac{\partial \mathbf L}{\partial \xx_i} \uu_i \Big)\vv,
\end{align}
which we refer to as the \textit{merged lower bound (of the eigenvalues' rate of change)} as it bounds the rate of change of all the eigenvalues $\{ \lambda_p \}_{p\in \mathcal I}$ at $\xx$ for a given~$\uu$.

We are interested in using the merged lower bound to replace the discontinuous function $\mu_2$ used in $\UU_\cm$, in order to avoid sudden changes in its value. For instance, noticing how the eigenspace $\VV_2$ expands into $\VV_{\until{2:3}}$ when $\lambda_2=\lambda_3$, we want to replace $\mu_2$ with $\mu_{\until{2:3}}$. This way, we avoid the abrupt change in the connectivity maintenance constraint map that occurs when $\lambda_2=\lambda_3$. However, with this approach, a discontinuity would still arise when $\lambda_3=\lambda_4$ since the eigenspace of $\lambda_4$ is not considered in the merged eigenspace. To address this, we can indeed use $\mu_{\until{2:n}}$, corresponding to the merged eigenspace of all nonzero eigenvalues, as stated in the following result. 

\begin{theorem}\label{thm:strict}
\longthmtitle{Strict Connectivity Constraint Map for Continuous Controller}
Consider the multi-robot system~\eqref{sys:full}. Given a locally Lipschitz extended class $\mathcal K$ function~$\alpha$, define the constraint map
\begin{equation}\label{eq:str_constraint}
\UU_{\str}(\xx) := \setdefb{\uu\in \real^N}{
\mu_{\until{2:n}}(\xx,\uu) \geq -\alpha(\lambda_2(\xx)-\varepsilon)} .
\end{equation}
If, for each $\xx$, there exists a control input $\uu\in\real^N$ that strictly satisfies the constraint map $\UU_{\str}\cap \UU_{\nom}(\xx)$, then under Assumptions~\ref{assump:continuity} and~\ref{assump:convexity}, the optimization-based controller 
\begin{equation}\label{eq:str_controller}
\kk_\str(\xx) := \underset{\uu\in \UU_\str\cap \UU_\nom(\xx)}{\arg\!\min}~J(\xx,\uu)
\end{equation}
is continuous on $\CC_\varepsilon$. In addition, the closed-loop feedback $\uu=\kk_\str(\xx)$ renders $\lambda_2(\xx(t))\geq \varepsilon$ at all time, ensuring that network connectivity is maintained, for any given initial condition $\xx_0\in\CC_\varepsilon$.~\hfill~\qed
\end{theorem}

While Theorem~\ref{thm:strict} provides a solution to Problem~\ref{prob:connect}, it is undoubtedly conservative. By design, the constraint map $\UU_{\str}$ bounds the rate of change of $\lambda_2$ as if it always has the highest possible multiplicity of~$n-1$ for a connected robot configuration. As a result, for situations when the multiplicity of $\lambda_2$ is unlikely to change, e.g., when $\lambda_2$ is far apart from $\lambda_3$, the design is conservative. This conservatism is also illustrated later in our simulations in Section~\ref{sec:experiment}.

To be less conservative, our next design takes into account how far the multiplicity of the eigenvalues is from changing. Instead of defining a NCBF constraint map for only $\lambda_2$, the design considers NCBFs for all the nonzero eigenvalues. We then replace each $\mu_m$ with the merged lower bound~$\mu_{\until{2:m}}$. Formally, for each $m\in \until{2:n}$, consider the constraint maps,
$$
\UU_{\until{2:m}}(\xx) := \setdefb{\uu\in \real^N}{
\mu_{\until{2:m}}(\xx,\uu) \geq -\alpha(\lambda_m(\xx)-\varepsilon)}
$$
with a locally Lipschitz extended class $\mathcal K$ function~$\alpha$ and a constant $\varepsilon\in\realpos$. The aggregations of the constraint maps of this form gives rise to our design for connectivity maintenance.

\begin{theorem}\longthmtitle{Aggregate Connectivity Constraint Map for Continuous Controller}\label{thm:main}
Consider the multi-robot system~\eqref{sys:full}. Given a locally Lipschitz extended class $\mathcal K$ function~$\alpha$, define the constraint map
\begin{equation}\label{eq:main_constraint}
\UU_{\agg}(\xx) := \bigcap_{m\in \until{2:n}}\UU_{\until{2:m}}(\xx).
\end{equation}
If, for each $\xx$, there exists a control input $\uu\in\real^N$ that strictly satisfies the constraint map~~$\UU_{\agg}\cap \UU_{\nom}(\xx)$, then under Assumptions~\ref{assump:continuity} and~\ref{assump:convexity}, the optimization-based controller
\begin{equation}\label{eq:main_controller}
\kk_\agg(\xx) := \underset{\uu\in \UU_\agg(\xx)\cap \UU_\nom(\xx)}{\arg\!\min}~J(\xx,\uu)
\end{equation}
is continuous on $\CC_\varepsilon$. In addition, the closed-loop feedback $\uu=\kk_\agg(\xx)$ renders $\lambda_2(\xx(t))\geq \varepsilon$ at all time, ensuring that network connectivity is maintained, for any given initial condition $\xx_0\in\CC_\varepsilon$.~\hfill~\qed
\end{theorem}

The idea behind the design of the aggregate constraint~\eqref{eq:main_constraint} is as follows. Consider a state $\xx$ where $\lambda_{m-1}(\xx) = \lambda_m(\xx)$. At this state, $\UU_{\until{2:m-1}}(\xx)$ abruptly shrinks to $\UU_{\until{2:m}}(\xx)$ due to the value of the merged lower bound~$\mu_\until{2:m-1}(\xx,\uu)$ dropping to that of $\mu_\until{2:m}(\xx,\uu)$, for any given $\uu$. Nevertheless, the constraint map $\UU_{\until{2:m}}$ is also considered in the aggregate constraint map~$\UU_\agg$, and the fact that it experiences no abrupt change there is enough to prevent $\UU_\agg$ from changing abruptly at that state.

Both constraint maps~\eqref{eq:str_constraint} and~\eqref{eq:main_constraint} ensure continuity of the corresponding optimization-based controller and solve Problem~\ref{prob:connect}. In general, for $m \in \until{2:n}$, one has $\UU_\str \subseteq \UU_{\until{2:m}}$ because $\mu_{\until{2:n}}\leq \mu_{\until{2:m}}$ and $\lambda_m \geq \lambda_2$. Therefore, $\UU_\str \subseteq \UU_\agg$, with equality holding on those states where $\lambda_n (\xx) = \lambda_2 (\xx)$. Consequently, $\UU_\agg$ imposes less conservative constraints than~$\UU_\str$. This is because the aggregate constraint~$\UU_\agg$ only gradually becomes stricter as the gap between each eigenvalue and the lowest $\lambda_m-\lambda_2$ gets smaller, unlike the strict constraint~$\UU_\str$ that is agnostic to the gap.

\begin{remark}\longthmtitle{Strictly Satisfying Feasible Controls Requirement}
\label{rmk:feasible}
{\rm We note that both Theorems~\ref{thm:strict} and~\ref{thm:main} require the existence, at each $\xx$, of a control $\uu$ strictly satisfying the corresponding constraint map. This is our conceptualization of the fact that, in order for Problem~\ref{prob:connect} to be solvable, there must exist at each state a control that can simultaneously maintain connectivity and satisfy the nominal constraints. The choice of class $\mathcal{K}$ function also provides flexibility in this regard because, if a control exists that satisfies the constraints at $\xx$ for $\alpha_1$, then the same control strictly satisfies the constraints for $\alpha_2$ with $\alpha_1 < \alpha_2$, as long as $\lambda_2(\xx)\neq \varepsilon$. Finally, as we show later in our analysis (cf. Lemma~\ref{lem:artstein}), the existence of strictly satisfying feasible control at each state is enough to guarantee the existence of a continuous controller. While this latter condition would be enough to establish Theorems~\ref{thm:strict} and~\ref{thm:main}, the existence of strictly satisfying feasible control is easier to check as it consists of a pointwise condition at each network state $\xx$, instead of the analysis across the states required to ensure continuity.~\hfill $\bullet$
}
\end{remark}

\begin{remark}\longthmtitle{Computation of Proposed Controllers}
\label{rmk:computation}
{\rm For each $\xx$, the computation of the controllers $\kk_\str$ and $\kk_\agg$ are convex optimization problems (as we show later, the constraint maps are convex-valued, and the cost function $J$ is convex by assumption). This means that one can utilize the wide variety of existing methods and solvers available for convex optimization, cf. \citep{SB-LV:09,RTR:70}, to compute the controllers. In implementing these methods, one must pay attention to the fact that obtaining the value of each merged lower bound function~$\mu_{\until{2:m}}$ is itself an optimization problem. Nevertheless, this can be addressed by casting the computation of the merged lower bounds as an eigenvalue problem. To see why this is so, note the following relationship
\begin{align*}
& \mu_{\until{2:m}}(\xx,\uu) = \min_{\vv\in\VV_{\until {2:m}}(\xx)}\vv^\top \Big(\sum_{i\in\until{N}} \frac{\partial \mathbf L}{\partial \xx_i} \uu_i \Big)\vv \\
& \qquad = \min_{\xi \in \SS^{m-1}}
\xi^\top [\vv]^\top_{2:m}(\xx) \Big(\sum_{i\in\until{N}} \frac{\partial \mathbf L}{\partial \xx_i} \uu_i \Big)[\vv]_{2:m}(\xx)\xi\\
& \qquad :=\min_{\xi \in \mathbb{S}^{m-1}}
\xi^\top  Z_m(\xx,\uu)\xi ,
\end{align*}
where $[\vv]_{2:m}(\xx)$ is the matrix created by concatenating orthonormal eigenvectors of $\{\lambda_p\}_{p\in\until{2:m}}$. It then follows that $\mu_{\until{2:m}}(\xx,\uu)$ is the minimum eigenvalue of the matrix $Z_m(\xx,\uu)$ defined above. This formulation as eigenvalue problem is advantageous for two reasons: it makes the evaluation of the function easy using standard linear algebraic routines and, for gradient-based optimization methods, it facilitates the computation of the generalized gradient of the merged lower bound.~\hfill $\bullet$
}
\end{remark}

\section{Technical Analysis of the Proposed Solutions}\label{sec:analyses}
This section provides the proofs of the results presented in Section~\ref{sec:results}. Before presenting them, we establish a number of auxiliary results that characterize the properties of the merged lower bounds involved in the construction of the constraint set-valued maps.

\subsection{Properties of Merged Lower Bounds}
We first examine the properties of functions $\mu_{\mathcal I}$ of the form~\eqref{eq:mu-I} defining our proposed constraint sets. The definition of such functions relies critically on the normalized merged eigenspace~$\VV_{\mathcal I}$. The following result characterizes the continuity properties of the latter.

\begin{theorem}
\label{thm:eig}
\longthmtitle{Continuity of Normalized Merged Eigen\-spaces}
Let $\map{\mathbf{L}}{\real^N}{\sym_n}$ be a continuous function. Given $\mathcal I\subset \until{n}$, the normalized merged eigenspace~$\VV_{\mathcal I}$ is continuous at any $\xx$ such that $\lambda_i(\xx)\neq \lambda_j(\xx)$ for all $i\in \mathcal I$ and $j\not\in\mathcal I$, i.e., where none of the considered eigenvalues is equal to any of the remaining eigenvalues.~\hfill$\qedA$
\end{theorem}

Due to its length, the proof of this result is provided in the Appendix. Building on this result, the continuity of the merged lower bounds follows from a direct application of the Berge Maximum Theorem~\citep[Thm.~17.31]{CDA-KCB:99}.

\begin{corollary}\label{cor:mu-I_cont}\longthmtitle{Continuity of Merged Lower Bounds}
Given $\mathcal I\subset \until{n}$, the function $\mu_{\mathcal I}$ is continuous at any $(\xx,\uu)$ such that  $\lambda_i(\xx)\neq \lambda_j(\xx)$ for all $i\in \mathcal I$ and $j\not\in\mathcal I$.~\hfill$\qedA$
\end{corollary}

In particular, we consider indices $\mathcal I = \until{2:m}$ of ordered eigenvalues on the domain where the graph remains connected $\CC_\varepsilon$ (i.e., where $\lambda_1(\xx)\neq \lambda_2(\xx)$). Thus, $\mu_{\until{2:m}}$ is continuous at any $\xx$ such that $\lambda_m(\xx) \neq \lambda_{m+1}(\xx)$, and $\mu_\until{2:n}$ is continuous everywhere on~$\CC_\varepsilon\times \real^N$. 

Besides continuity of $\mu_{\until{2:m}}$, another crucial property to show is convexity of the constraint maps $\UU_\str$ and $\UU_\agg$. To this end, we establish the concavity property of the merged lower bounds.

\begin{lemma}\longthmtitle{Concavity of Merged Lower Bounds}\label{lem:constraint_convexity}
For any $\mathcal I \subseteq \until n$, $\mu_{\mathcal I}$ is concave in $\uu$. Consequently, the constraint maps $\UU_\str$ and $\UU_\agg$ are convex-valued. 
\end{lemma}
\begin{pf}
Given any $\uu^1,\uu^2\in\real^N$ and $0\leq \gamma\leq 1$, we have
\begin{align*}
\mu_{\mathcal I}(\xx,&\gamma \uu^1 + (1-\gamma) \uu^2) \\
&= \min_{\vv\in\VV_{\mathcal I}(x)}\vv^\top \Big(\sum_{i\in\until{N}} \frac{\partial \mathbf L}{\partial \xx_i} (\gamma \uu^1_i+(1-\gamma) \uu^2_i) \Big)\vv 
\\
&\geq\min_{\vv\in\VV_{\mathcal I}(\xx)}\Big(\gamma \vv^\top \Big(\sum_{i\in\until{N}} \frac{\partial \mathbf L}{\partial \xx_i} \uu^1_i \Big)\vv \Big)
\\
&\quad +\min_{\vv\in\VV_{\mathcal I}(\xx)}\Big( (1-\gamma) \vv^\top \Big(\sum_{i\in\until{N}} \frac{\partial \mathbf L}{\partial \xx_i} \uu^2_i \Big)\vv  \Big)
\\
&= \gamma \mu_{\mathcal I}(\xx,\uu^1) +(1-\gamma) \mu_{\mathcal I}(\xx,\uu^2).
\end{align*}
Therefore, $\mu_{\mathcal I}$ is concave in $\uu$.~\hfill~\qed
\end{pf}

Having established the continuity and concavity properties of the merged lower bounds $\mu_{\mathcal I}$, we next turn our attention to characterize the properties of the constraint maps.

\begin{remark}\longthmtitle{More General System Dynamics}
{\rm Note that Lemma~\ref{lem:constraint_convexity} is the only instance in our technical treatment where we have exploited the specific form of the single-integrator dynamics~\eqref{sys:full}. In fact, all of our results are valid for more general system dynamics, as long as one can prove the concavity of the merged lower bound. This is the case, for instance, for control-affine systems, for which our results hold. We have kept the discussion limited to single-integrator dynamics for simplicity of presentation.
~\hfill $\bullet$
}
\end{remark}

\subsection{Equivalent Constraint Maps}
In general, the constraint maps $\UU_\str$ and $\UU_\agg$ might not be UHC because they are unbounded. To make sure the requirements of Lemma~\ref{lem:UHC_req} as well as Berge Maximum Theorem are met, we explain here how to consider, following Section~\ref{sec:problem_statement}, equivalent constraint maps that are compact-valued. This procedure involves using sublevel sets of the cost function~$J$, which are compact due to Assumption~\ref{assump:convexity}. In order to do so, we require a feasible control function $\xx \mapsto \bar \uu(\xx)$ to define $\JJ_{\bar \uu}$ as in~\eqref{eq:J_sublvl}. Note, importantly for our purposes, that the function $\bar \uu$ must be continuous so that $\JJ_{\bar \uu}$ is also continuous. The next result shows that, under the assumptions of Theorems~\ref{thm:strict} and~\ref{thm:main}, such continuous feasible control function always exists.

\begin{lemma}\longthmtitle{Generalization of Artstein's Theorem}\label{lem:artstein}
Consider a set-valued map $\setmap{\UU}{\real^N}{\real^M}$ defined with a vector-valued function $\map{g}{\real^N}{\real^M}$ as
$$
\UU(\xx) = \setdef{\uu\in\real^M}{g(\xx,\uu)\leq 0}.
$$
If $g$ is continuous and $\UU$ is convex-valued, and, for each $\xx$, there exists a control input $\uu$ that strictly satisfies $\UU(\xx)$, then there exists a $\mathcal C^\infty$ function $\map{\bar \uu}{\real^N}{\real^M}$ such that $\bar \uu(\xx)\in\UU(\xx)$.
\end{lemma}
\begin{pf}
For each $\xx$, let $\uu_\text{int}(\xx)$ denote the control input such that
$g(\xx,\uu_\text{int}(\xx)) < 0$. Due to continuity of $g$, there exists a neighborhood of $\xx$, denoted by $\mathcal W(\xx)$, such that $\uu_\text{int}(\xx)\in\UU(\xx')$ for all $\xx'\in \mathcal W(\xx)$. The collection of $\{\mathcal W(\xx)\}_{\xx\in\real^N}$ is an open cover for $\real^N$. Then, because we deal with a Euclidean space that is a differentiable manifold, there exists a countable partition of unity $\{\psi_j\}$ subordinate to the cover, cf. \citep[Thm. 1.11]{FWW:89}. In other words, for each~$j$, there exists an~$\xx$ such that $\text{supp}(\psi_j)$ is a subset of $\mathcal W(\xx)$, each of which has an associated control $\uu_\text{int}^j\in \UU(\xx)$ for $\xx\in \text{supp}(\psi_j)$. Then we define $\bar \uu(\xx)=\sum_j\psi_j(\xx)\uu_\text{int}^j$, which satisfies the statement due to convexity of the map~$\UU$.~\hfill~\qed
\end{pf}

Lemma~\ref{lem:artstein} is a generalization of Artstein's Theorem~\citep[Thm. 4.1]{ZA:83} on the existence of a continuous controller given a control Lyapunov function. The proof of the result, included here for completeness, is also a slight modification of the original proof. Because the functions defining $\UU = \UU_\str$ are continuous, we can directly apply Lemma~\ref{lem:artstein}. On the other hand, $\UU = \UU_\agg$ is defined with discontinuous functions; nevertheless, from its construction, one can still employ the argument presented in the proof of Lemma~\ref{lem:artstein} (i.e., there exists $\uu_\text{int}$ at each $\xx$ belonging to $\UU_\agg(\xx')$ for all $\xx'$ in a neighborhood $\mathcal W$ of $\xx$, and so on). As a result, for each of the cases  $\UU = \UU_\str$ and $\UU = \UU_\agg$, there exists a continuous feasible control function $\bar \uu$, which we use to define the corresponding set-valued map~$\JJ_{\bar \uu}$. This map is convex-valued and compact-valued due to it being a sublevel set of a strictly convex function $J$, cf. Assumption~\ref{assump:convexity}. Then according to Lemmas~\ref{lem:UHC_req} and~\ref{lem:LHC_req}, it is also continuous due to the functions $\bar \uu$ and $J$ being continuous, cf. Assumption~\ref{assump:continuity}. We then consider the intersections $\UU_\str \cap \UU_\nom\cap\JJ_{\bar \uu}$ and~$\UU_\agg\cap \UU_\nom\cap \JJ_{\bar \uu}$, where the inclusion of  $\JJ_{\bar \uu}$ make these constraint maps compact-valued. For the purpose of our analysis, we equivalently define $\kk_\str$ and $\kk_\agg$ with these constraint maps as the constraint to the optimization.

\subsection{Continuity of the Connectivity Maintenance Controllers}

With the preparations from prior sections, we are now ready to prove our results on continuity of $\kk_\str$ and $\kk_\agg$.

\begin{pf}[Proof of Theorem~\ref{thm:strict}]
Consider the constraint set $\UU_\str \cap \UU_\nom \cap \JJ_{\bar \uu}$. We note the following: (i) all the functions defining the constraint map are continuous due to Assumption~\ref{assump:continuity} and $\mu_\until{2:n}$ being continuous everywhere (on $\CC_\varepsilon\times\real^N$); (ii) the map is convex-valued because all intersecting maps are convex-valued; (iii) the map has a nonempty interior by assumption; (iv) the map is compact-valued because the intersecting maps are closed-valued and $\JJ_\uu$ is compact-valued. Thus, we may apply Lemmas~\ref{lem:UHC_req} and~\ref{lem:LHC_req}, to show continuity of this constraint map. By Berge Maximum Theorem~\citep[Thm. 17.31]{CDA-KCB:99}, $\kk_\str$ is a continuous function as stated. Lastly, from the relationship 
$$
\min \Lie_{f_{\sint}}\lambda_2(\xx,\uu) = \mu_2(\xx,\uu) \geq \mu_\until{2:n}(\xx,\uu),
$$
it follows that $\kk_\str(\xx)\in \UU_\str(\xx)\subseteq \UU_\cm(\xx)$. As a result, Lemma~\ref{lem:connect} ensures $\lambda_2(\xx(t))\geq \varepsilon$, and the proof concludes.~\hfill~\qed
\end{pf}
We next prove the continuity result for $\kk_\agg$, which is more complicated due to the merged lower bounds used not being continuous everywhere.
\begin{pf}[Proof of Theorem~\ref{thm:main}] Consider the constraint map $\UU_\agg\cap\UU_\nom\cap\JJ_{\bar u}$. Because each $\mu_\until{2:m}$ is not continuous everywhere, we can only conclude continuity using Lemmas~\ref{lem:UHC_req} and~\ref{lem:LHC_req} wherever $\mu_\until{2:m}$ are continuous for all $m\in\until{2:n}$. For the remaining states, we show continuity of the constraint map by proving separately below that it is UHC and LHC. Note that once we prove continuity, the theorem statements are established analogously as we did in the proof of Theorem~\ref{thm:strict}

\emph{Upper Hemicontinuity:} We begin by consider the partial constraint map $\JJ_{\bar \uu} \cap\UU_\until{2:n}$. This set-valued map is continuous on~$\CC_\varepsilon$ because of the continuity of $\mu_\until{2:n}$. Consider its intersection with $\JJ_{\bar \uu}\cap \UU_\until{2:n-1}$. At the states where $\lambda_n(\xx)=\lambda_{n-1}(\xx)$, notice that $\UU_\until{2:n}(\xx) = \UU_\until{2:n-1}(\xx)$, so the intersection  $\JJ_{\bar \uu}\cap (\bigcap_{m\in\until{n-1:n}}\UU_\until{2:m})(\xx)$ is exactly the same set as $\JJ_{\bar \uu} \cap\UU_\until{2:n}(\xx)$ at those $\xx$. For other states~$\xx$, we know that the former map is a subset of the latter. Then, directly from the definition of UHC for $\JJ_{\bar \uu} \cap\UU_\until{2:n}$, we can conclude UHC for the intersection $\JJ_{\bar \uu}\cap (\bigcap_{m\in\until{n-1:n}}\UU_\until{2:m})$ at $\xx$ where $\lambda_n(\xx)=\lambda_{n-1}(\xx)$. 
Elsewhere, the intersection can be proven UHC directly via Lemma~\ref{lem:UHC_req}, so it is continuous everywhere on~$\CC_\varepsilon$. With the same line of reasoning, we can continue to show by induction that $\JJ_{\bar u} \cap \UU_\agg$ is UHC on~$\CC_\varepsilon$. Then intersecting with $\UU_\nom$, we conclude the set-valued map~$\UU_\agg\cap\UU_\nom\cap\JJ_{\bar u}$ is UHC from Lemma~\ref{lem:insct-uhc}.

\emph{Lower Hemicontinuity:} We begin by defining the following auxiliary set-valued maps for  $m\in\until{2:n}$,
$$
\HH_m(\xx) = \setdefb{\uu\in \real^n}{\mu_{\until{2:m}}(\xx,\uu)\geq -\alpha(\lambda_{m-1}(\xx)-\varepsilon)}.
$$
By definition, $\HH_m(\xx)\subseteq \UU_\until{2:m-1}(\xx) $ because $\mu_\until{2:m}(\xx,\uu) \leq \mu_\until{2:m-1}(\xx,\uu)$, and $\HH_m(\xx) \subseteq \UU_\until{2:m}(\xx)$ because $\lambda_{m-1}(\xx) \leq \lambda_{m}(\xx)$. In addition, note that $\HH_m$ is convex-valued because the merged lower bound~$\mu_{\until{2:m}}$ in concave in $\uu$, cf. Lemma~\ref{lem:constraint_convexity}, and it has a 
nonempty interior as it is a subset of $\UU_\until{2:m}$, which has a nonempty interior by assumption. Then, by Lemma~\ref{lem:LHC_req} it is LHC for all $\xx\in\CC_\varepsilon$ where $\lambda_m(\xx) \neq \lambda_{m+1}(\xx)$ (with $\HH_n$  continuous everywhere on $\CC_\varepsilon$).

We prove LHC of $\UU_\agg$ by induction. We start by considering the maps $\UU_\until{2:n}$ and $\HH_n$, which are both LHC on~$\CC_\varepsilon$. We then consider the intersection $\UU_\until{2:n}$ with  $\UU_\until{2:n-1}$. For $\xx$ where $\lambda_n(\xx)=\lambda_{n-1}(\xx)$, the two eigenvalues share the same eigenspaces. Thus, it is also the case that $\mu_\until{2:n}(\xx)=\mu_\until{2:n-1}(\xx)$, and we find that $\HH_n(\xx) = \UU_\until{2:n}(\xx)\cap \UU_\until{2:n-1}(\xx)$ for all $\xx$ where the two eigenvalues are equal. From this and the fact that $\HH_n$ is a subset of $\UU_\until{2:n}\cap \UU_\until{2:n-1}$ in general, we can use the LHC of $\HH_n$, at $\xx$ where $\lambda_n(\xx)=\lambda_{n-1}(\xx)$ to deduce LHC for  $\UU_\until{2:n}\cap \UU_\until{2:n-1}$ there. Elsewhere, the set $\UU_\until{2:n}\cap \UU_\until{2:n-1}$ can be proven continuous directly from Lemma~\ref{lem:LHC_req}, so it is LHC everywhere on $\CC_\varepsilon$. Then using Lemma~\ref{lem:insct-lhc}, we also deduce that the intersection $\HH_n \cap (\UU_\until{2:n}\cap \UU_\until{2:n-1})$ is LHC on $\CC_\varepsilon$. 

To continue with the induction proof, assume the set-valued maps 
$$
\bigcap_{m\leq p\leq n} \UU_{[2:p]}~\text{and}~\HH_{m} \cap \bigcap_{m\leq p\leq n} \UU_{[2:p]}
$$
are LHC. Then we can follow the arguments above to also deduce that their intersections with $\UU_\until{2:m-1}$ are also LHC. Hence, $\UU_\agg$ is LHC. Then the LHC of the intersection $\UU_\agg \cap \UU_\nom \cap \JJ_{\bar \uu}$ follows via Lemma~\ref{lem:insct-lhc}, concluding the proof.~\hfill~\qed
\end{pf}

\section{Simulations and Experimental Validation}\label{sec:experiment}
In this section we report the simulations and the experiment we have carried out to verify the effectiveness of the proposed controller. We consider a resource gathering problem with a group of four ($n=4$) robots, moving in a two-dimensional space ($d_r=2$ for all agents). Each robot is tasked with visiting its own target region. If the robots were to individually move directly to their targets, the network will be disconnected. Therefore, we prioritize the order in which the robots reach their targets and use our proposed controller to maintain the connectivity among them. We consider the mission accomplished when the target location is visited by the corresponding robot, and we change the task prioritization to the next robot.

The nominal controller carries each robot towards the corresponding target with a conical potential field:
\begin{equation}
    \label{eq:our_u}
    u_{\nom,r}(x_r) = v_{\nom} \frac{e_r(x_r)}{\|e_r(x_r)\|},~\forall r\in \until{n},
\end{equation}
where $x_r$ is the position of the $r$-th robot and $e_r(x_r) = x_{\target,r} - x_r$ is the error between the center of the robot's target region and its position, and $v_{\nom} \in \realpos$ is a constant velocity parameter. By denoting $\kk_\nom(\xx) = [u_{\nom,1}(x_1)^\top,\dots,u_{\nom,4}(x_4)^\top]^\top$, our cost function
\begin{equation}
    \label{eq:our_J}
    J(\xx,\uu) = \|\uu - \kk_{\nom}(\xx) \|^2
\end{equation}
measures the deviation of the control decision from the nominal controller. In order to ensure that our prioritized robot, indexed $P$, makes progress towards its target, we enforce the following constraint map,
\begin{equation}
    \label{eq:our_U}
    \UU_{\nom}(\xx) = \setdefb{\uu \in \real^N}{k v_{\nom}^2 - u_{\nom,P}^{\top} u_P \leq 0},
\end{equation}
where $k \in \realpos$ is a constant parameter to restrict how much $u_P$ should point in the direction of~$u_{\nom,P}$. Once the robot reaches its assigned target region, its $u_{\nom,P}$ is set equal to zero. This represents the fact that after having accomplished its task, the robot is relieved from its mission and prefers to conserve energy by not moving. Note that it continues to collaborate at maintaining the connectivity. Also after the prioritized robot achieves its mission, we change the value of $P$ to correspond to the index of the next robot that has yet to achieve its goal, whose task we want to prioritize. 

Note that the objective~\eqref{eq:our_J} and the nominal constraint~\eqref{eq:our_U} verify both Assumptions~\ref{assump:continuity} and~\ref{assump:convexity} for any prioritized robot (we disregard the jumps in  $\UU_{\nom}$ due to 
the transitions when a robot reaches its target region and the identity of the prioritized robot changes). In addition, we can verify strict feasibility of the constraints since constraint~\eqref{eq:our_U} only affects the prioritized robot $P$, and we can have all other robots move towards the prioritized robot to maintain connectivity (because we have assumed the single-integrator dynamics~\eqref{sys:full}). We show that our proposed controllers from Theorem~\ref{thm:strict} and~\ref{thm:main} are continuous for the duration between events when the prioritized robots achieve its goal. For both our simulations and our experiment, we use the projected saddle-point dynamics~\citep{AC-BG-JC:17} to solve the convex optimization problems at each given state $\xx$ along the trajectory. The saddle-point dynamics involve performing a gradient descent in the primal variable~$\uu$, and a gradient ascent in the dual variables (projected onto the positive orthant to ensure they remain nonnegative) associated with the constraints. We chose this gradient-based method because we know how to compute the gradient of the merged lower bounds that appear in our constraints, cf. Remark~\ref{rmk:computation}. The computation of our controllers is done in MATLAB\textsuperscript \textregistered.

\subsection{Simulations}\label{sec:simulations}
\begin{figure}[tb]
	\centering
	\begin{subfloat}[][\centering Aggregate controller]{
	\includegraphics[width=\columnwidth,keepaspectratio]{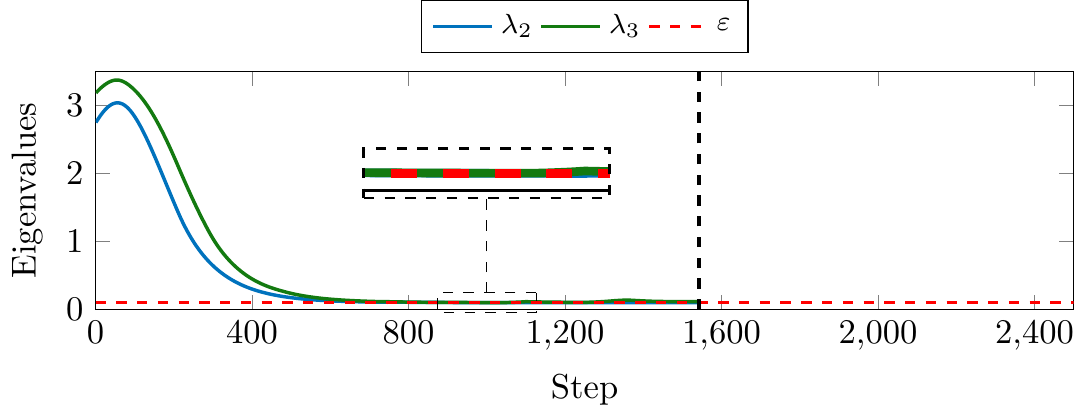}
    \label{fig:simulation_eigenvalues_1}}
	\end{subfloat} \\
    \begin{subfloat}[][\centering Strict controller]{
    \includegraphics[width=\columnwidth,keepaspectratio]{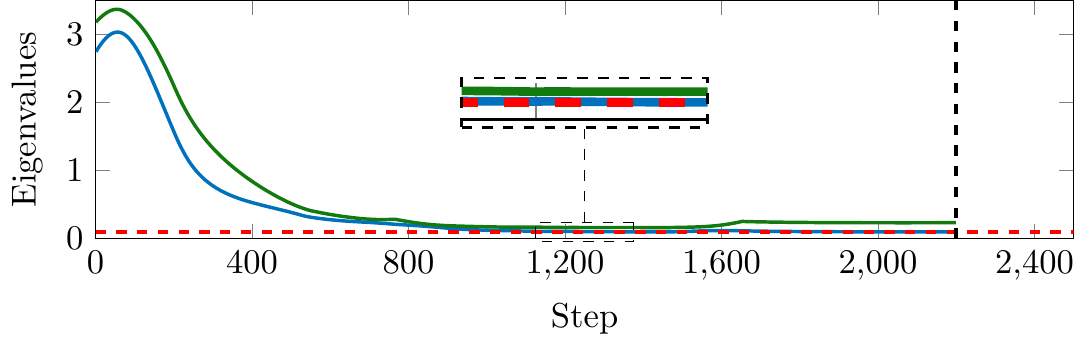}
    \label{fig:simulation_eigenvalues_2}}
	\end{subfloat} \\
	\begin{subfloat}[][\centering Discontinuous  controller]{
	\includegraphics[width=\columnwidth,keepaspectratio]{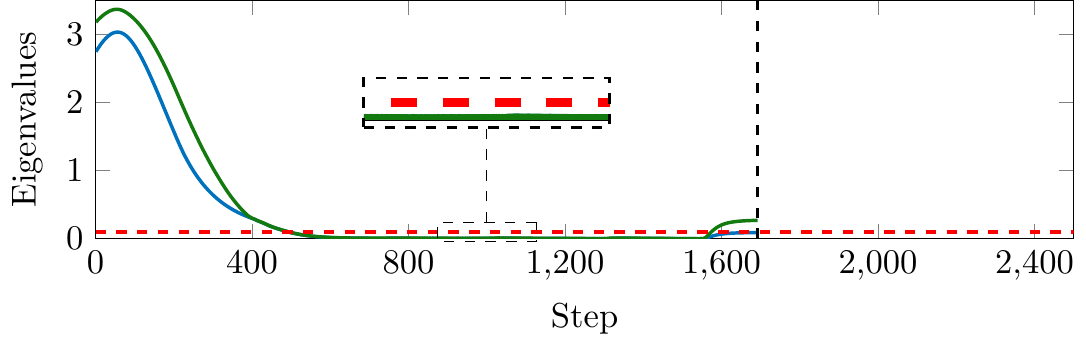}
    \label{fig:simulation_eigenvalues_3}}
	\end{subfloat}
	\caption{Eigenvalue evolution during the simulations under the different controllers. The \textit{dashed black} lines represent the end of the network task.}
	\label{fig:simulation_eigenvalues}
\end{figure}

Our simulations highlight the differences among the different controllers: $\kk_{\dis}$, defined in~\eqref{eq:dis_controller}, $\kk_\str$, defined in~\eqref{eq:str_controller}, and $\kk_\agg$, defined in~\eqref{eq:main_controller}, each with $\alpha(\lambda)=\lambda$. The initial positions, the robots' targets, and the parameters ($v_\nom = 0.5$, $k = 0.75$, $\varepsilon = 0.1$) are the same in each simulation. Fig.~\ref{fig:simulation_eigenvalues} reports the eigenvalues of the Laplacian matrix during the simulations. It is clear how both the aggregate (Fig.~\ref{fig:simulation_eigenvalues_1}) and the strict controller (Fig.~\ref{fig:simulation_eigenvalues_2})  maintain the connectivity constraint, unlike the discontinuous one that leads to disconnection (Fig.~\ref{fig:simulation_eigenvalues_3}). 

\begin{figure}[tb]
	\centering
	\begin{subfloat}[][\centering Aggregate controller]{
	\includegraphics[width=\columnwidth,keepaspectratio]{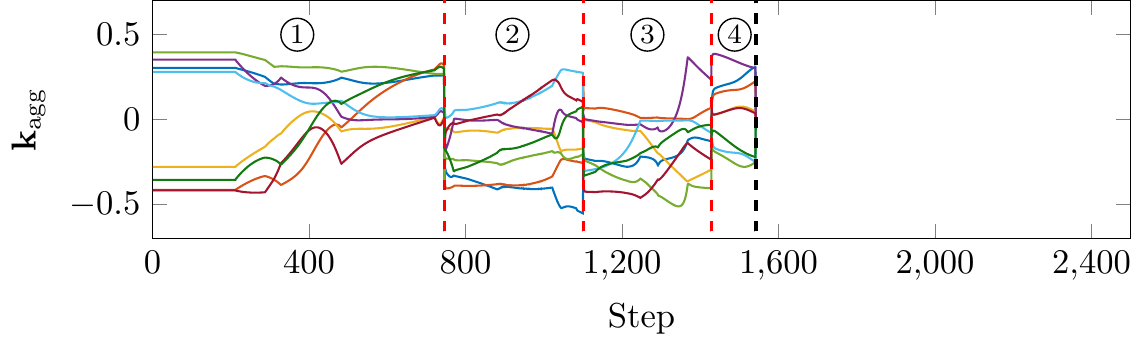}
	\label{fig:simulation_input_1}}
	\end{subfloat}\\
	\begin{subfloat}[][\centering Strict controller]{
	\includegraphics[width=\columnwidth,keepaspectratio]{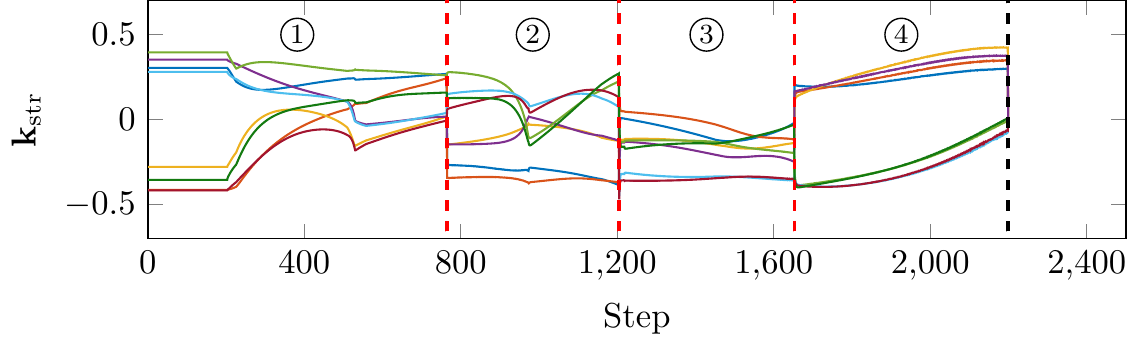}
	\label{fig:simulation_input_2}}
	\end{subfloat}\\
	\begin{subfloat}[][\centering Discontinuous controller]{
	\includegraphics[width=\columnwidth,keepaspectratio]{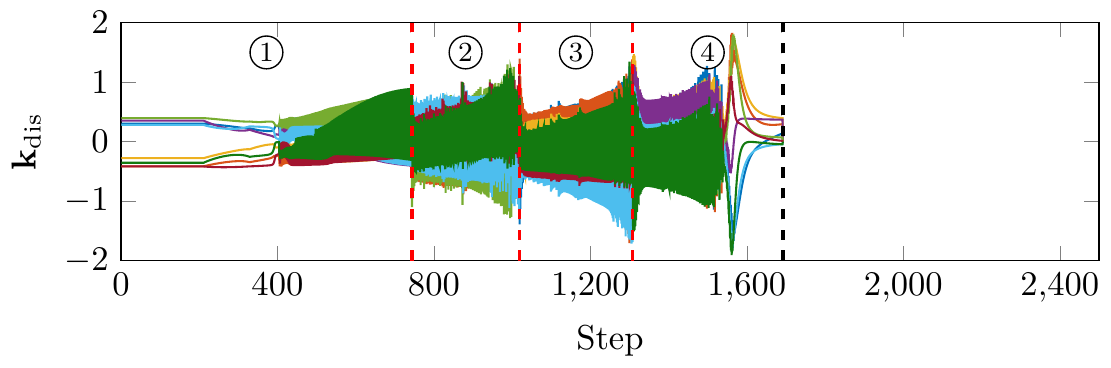}
	\label{fig:simulation_input_3}}
	\end{subfloat}
	\caption{Control inputs during the simulations of different controllers. We report all the components of the control input for each robot.}
	\label{fig:simulation_input}
\end{figure}

The overall performance corroborates our hypothesis that the aggregate controller will outperform the strict one. Figs.~\ref{fig:simulation_input_1} and~\ref{fig:simulation_input_2} show the continuous input produced by the aggregate and the strict controllers, and Fig.~\ref{fig:simulation_input_3} shows the discontinuous one generated by the discontinuous controller. It takes the aggregate controller (1542 steps) shorter time than the strict controller (2199 steps) to complete all the tasks in the simulation. We expected this result because the conservatism in the strict controller leads to fewer eligible robot formations. On the other hand, the aggregate controller allows the robots to get in better positions for the subsequent tasks. We note that there is no guarantee that the aggregate controller will always do better. This is because the two controllers generate two different trajectories and the strict controller may take the robots to better formations despite having fewer options available. Nevertheless, we believe that the aggregate controller will outperform the stricter one in general, based on our earlier reasoning. Regarding the discontinuous controller, we observe chattering throughout most of the trajectory. Such behavior is undesirable because it is difficult to implement in real systems, and it may also cause other unexpected issues. For example, despite theoretically being the least conservative out of the three controllers, the discontinuous controller does not complete the tasks before the aggregate one.  

\begin{figure}
    \centering
    \begin{subfloat}[][\centering Aggregate controller]{
    \includegraphics[width=\columnwidth,keepaspectratio]{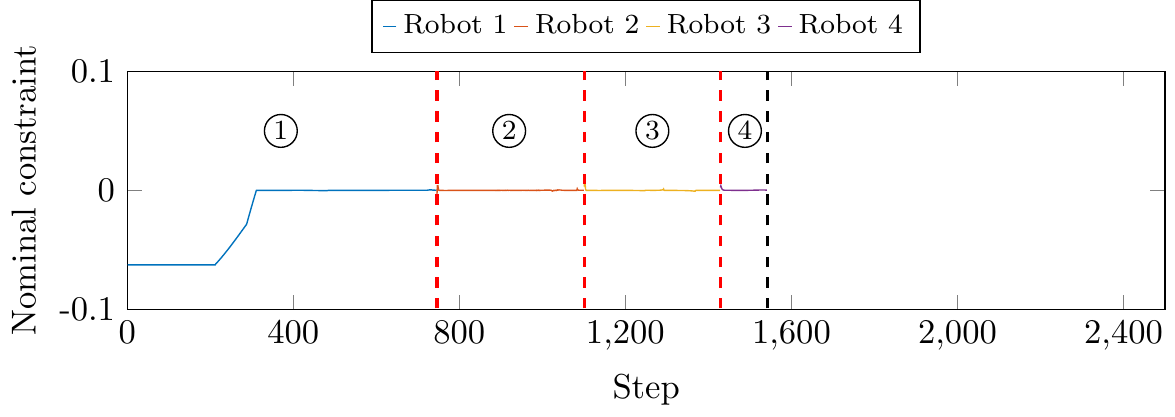}
        \label{fig:simulation_priority_1}}
    \end{subfloat}\\
    \begin{subfloat}[][\centering Strict controller]{
    \includegraphics[width=\columnwidth,keepaspectratio]{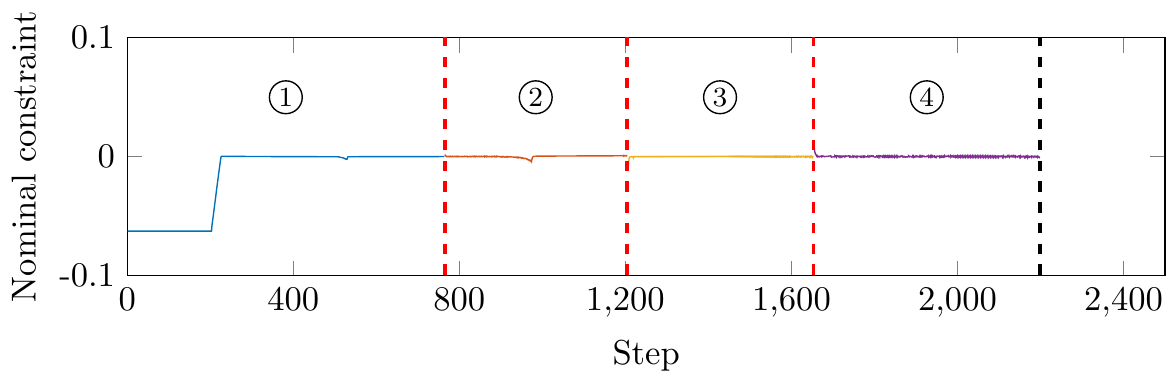}
        \label{fig:simulation_priority_2}}
    \end{subfloat}
    \caption{Nominal constraint~\eqref{eq:our_U} during the simulations. The \textit{dashed red} lines represent the instant in which the robot priority changes, due to the fact that a target has been reached. The number for each time slot corresponds to the robot with the highest priority.}
    \label{fig:simulation_priority}
\end{figure}

Next, Fig.~\ref{fig:simulation_priority} reports the evolution of the function defining the nominal constraint map~$\UU_\nom$ under the aggregate and the strict controllers. In the corresponding slot of time, the robot that has the target with the highest priority respects the constraint, while the others cooperate to maintain connectivity, minimally changing their nominal control law. We do not report the plot for the discontinuous controller as it is highly jittering, confirming what is already displayed in Fig.~\ref{fig:simulation_input_3}.

\subsection{Experimental Validation}

\begin{figure}[tb]
    \centering
    \includegraphics[width=0.8\columnwidth]{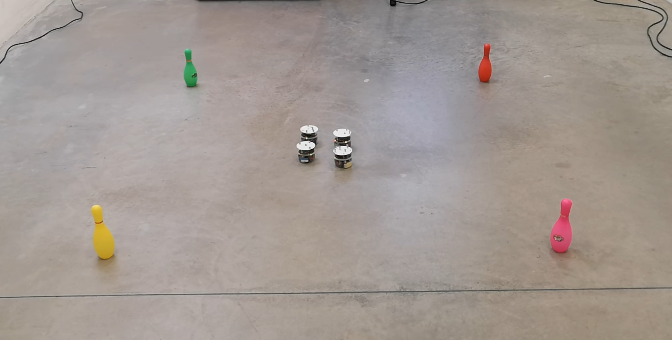}
    \caption{Experimental setup: 4 ePucks and their corresponding targets.}
    \label{fig:real_setup}
\end{figure}

We also carry out an experiment for the same resource gathering problem, cf. Fig.~\ref{fig:real_setup}. We use four small wheeled robots (ePucks) that are controlled via Bluetooth from a central unit that performs the calculations. The central unit is also connected to an Optitrack system, which provides the position of the robots. In order to transform the input calculated for the single-integrator dynamics to the unicycle dynamics of the robots, we use a simple input-output state-feedback linearization~\citep{GO-ADL-MV:02}. We tested only the proposed controller $\kk_{\agg}$, as the simulations  in Section~\ref{sec:simulations} verified that it is the best both in terms of performance and connectivity maintenance. We set the main parameters as $v_\nom = 0.1$, $k = 0.75$, and $\varepsilon = 0.3$. We report an example of the experiments in an accompanying video\footnote{https://youtu.be/auifIv-JR7o}.

\begin{figure}
    \centering
    \includegraphics[width=\columnwidth,keepaspectratio]{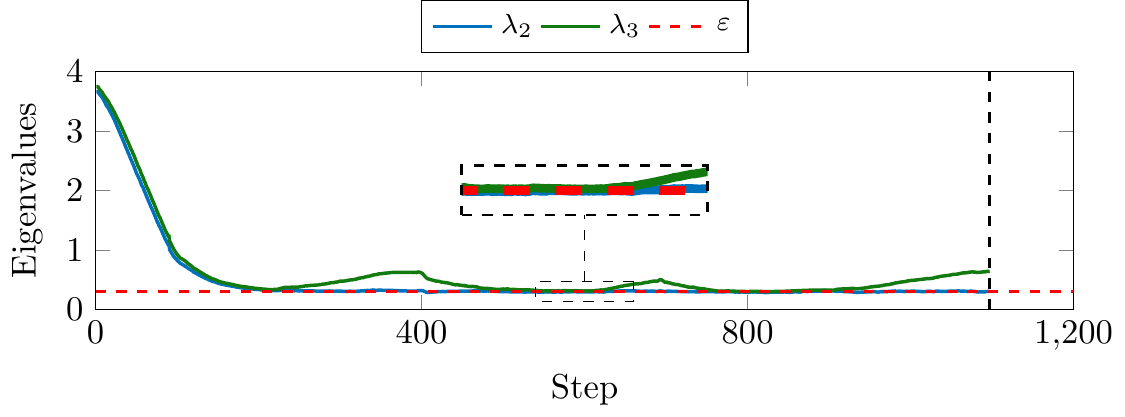}
    \caption{Eigenvalue evolution during the experiment. The \textit{dashed black} lines represents the end of the task.}
    \label{fig:real_eigenvalue}
\end{figure}

\begin{figure}[tb]
	\centering
	\includegraphics[width=\columnwidth,keepaspectratio]{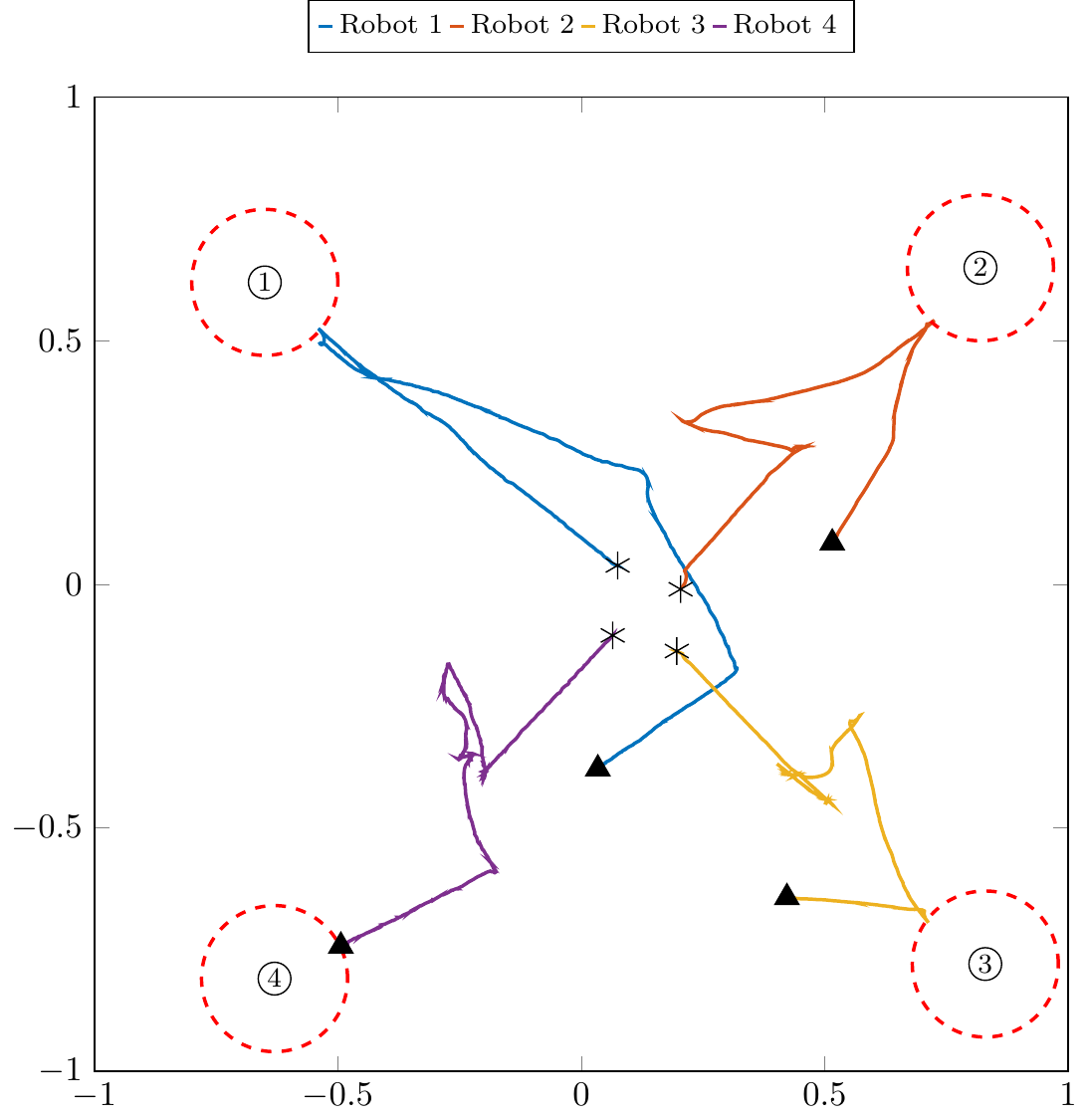}
	\caption{Trajectories followed by the robots during the experiment. The \textit{dotted red} circles represent the region where we consider the target reached (circle of 15 cm of radius around the target). The numbers represent the order of priority of the targets. The initial and final positions are reported with \textit{asterisks} and \textit{triangles}, respectively.}
	\label{fig:real_trajectory}
\end{figure}

\begin{figure}
    \centering
    \includegraphics[width=\columnwidth,keepaspectratio]{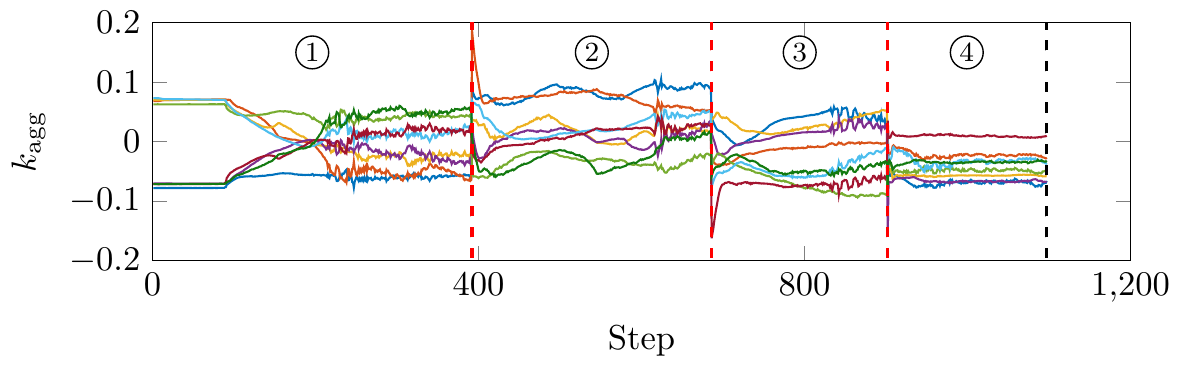}
    \caption{Control input applied to the robots in the experiment. We report all the components of the control input for each robot. These inputs have been transformed via input-output state-feedback linearization to be executed by the robots.}
    \label{fig:real_input}
\end{figure}

\begin{figure}
    \centering
    \includegraphics[width=\columnwidth,keepaspectratio]{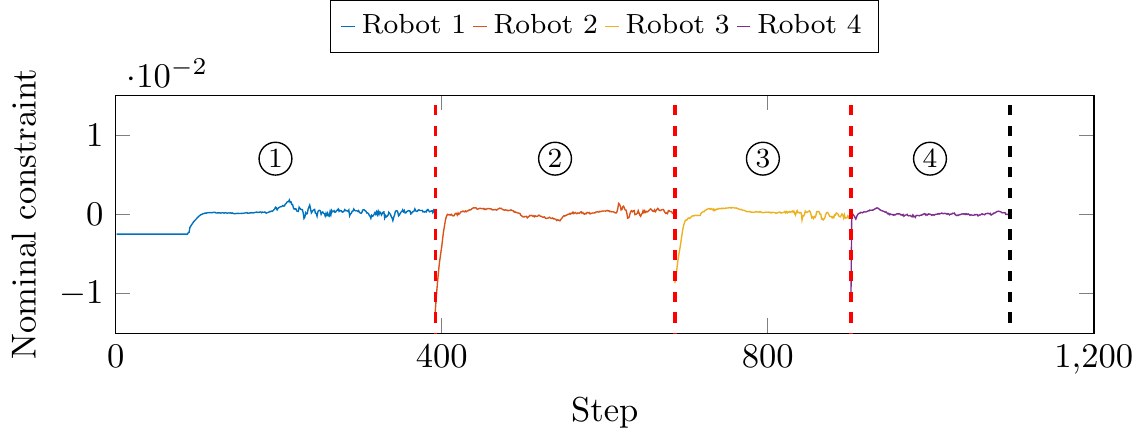}
    \caption{Nominal constraint~\eqref{eq:our_U} during the experiment. The \textit{dashed red} lines represent the instant in which the robot priority changes, due to the fact that a target has been reached. The number for each time slot corresponds to the robot with the highest priority.}
    \label{fig:real_priority}
\end{figure}

Fig.~\ref{fig:real_eigenvalue} reports the eigenvalue evolution during the experiment, further confirming the effectiveness of the proposed method in maintaining connectivity. Fig.~\ref{fig:real_trajectory} shows the trajectories followed by the robots, accomplishing the gathering task. It is evident how each target had been reached by the corresponding robot, and in the final positions (reported with \textit{triangles}) it is possible to see how the robots that have already reached their target cooperate to connectivity. Fig.~\ref{fig:real_input} shows the applied control inputs: here, the jittering is due both to the non-idealities introduced while using wheeled robots, which hardly instantaneously follow an omnidirectional dynamics, and the time needed for the calculation, which sometimes introduces a small delay. In fact, the time required to let the saddle-point dynamics converge is longer than the time needed to update the control input of the robots, which run at 10 Hz. For a given measurement of the state $\xx$, we run the saddle-point dynamics as many steps as possible so that the control $\uu$ solution is near optimal. Nevertheless, we report the computation of the gradient of the merged lower bounds as the bottleneck in the execution of the saddle-point dynamics. These are equivalent to solving eigenvalue problems, so they scale up cubically with the number of agents, $\mathcal O(N^3)$. In our future work, we plan to investigate ways of avoiding the calculation of these gradients when it is not necessary. Despite the limitations of the calculation and of the input of the robots, we achieve good performance also in satisfying the nominal constraint, cf. Fig.~\ref{fig:real_priority}.

\section{Conclusions}
We have considered the problem of maintaining network connectivity in multi-robot systems while satisfying nominal requirements that encode desired control objectives.  Our solution employs the algebraic connectivity of the interconnection topology as a nonsmooth control barrier function to produce additional constraints for the optimization-based synthesis of the controller that guarantee it is continuous and maintains network connectivity.  The technical approach fully embraces the nonsmooth nature of the algebraic connectivity and other spectral functions of the Laplacian matrix corresponding to the interconnection graph. This has led us to define two different continuous set-valued constraint maps, one that reasons with the merged lower bound of all the eigenvalues' rate of change at once and another, less conservative, that instead reasons over merged lower bounds of an increasing number of eigenvalues' rate of change.  We have illustrated the effectiveness of our approach in both simulation and experiment in a resource gathering multi-robot scenario. Future work will investigate the extension of our approach to ensure Lipschitzness of the controller, the application of the methodology proposed here to the synthesis of distributed controllers for connectivity maintenance,
the development of formal guarantees for robustness under perturbations and discrete-time implementations of our design, resource-aware design of aperiodic sample-and-hold implementations of the proposed controllers that do not require solving convex optimization problems at every state, and the extension of our results to more general nonlinear systems.

\appendix

\section{Appendix}
Here we provide several results on the continuity properties of eigenspaces, with the ultimate goal of establishing that the merged eigenspaces are continuous functions of the state as long as their dimensions remain constant, cf. Theorem~\ref{thm:eig}\footnote{Although this result is seemingly intuitive, we have not found it in the literature. There are results (e.g., \citep[Ch. 2.5.3]{TK:76}) that study the continuity properties of eigenvectors when their eigenvalues have multiplicity of one, a case where the eigenvectors can be viewed as a single-valued function. Instead, we investigate eigenspaces of eigenvalues with higher multiplicity, which requires set-valued analysis.}.

Given indices $\mathcal I \in \until{n}$, consider the merged eigenspaces $\VV_{\mathcal I}$.
For the purpose of analysis, instead of writing $\VV_\mathcal{I}$ as a span, we write out the full set definition as follows,
\begin{multline}\label{eq:eigspan}
  \VV_{\mathcal I}(\xx) = \setdefB{\vv\in \real^n}{
  (\mathbf L(x)-\lambda_i(x)\mathbf I)\xi_i = 0,~\forall i\in \mathcal I, \\
  \vv = \sum_{i\in \mathcal I} \cc_i\xi_i,~\cc\in\real^{\vert \mathcal I\vert},~\|\vv\|=1 }.
\end{multline}
For this set-valued map, we will show UHC and LHC separately.

\subsection{Upper Hemicontinuity of Merged Eigenspaces}
For the analysis of~\eqref{eq:eigspan}, it is convenient to use the eigenbasis as the coordinate system. Given a state $\xx^*\in \real^N$ at which we seek to prove continuity, let the matrix $\TT\in \real^{n\times n}$ be an orthonormal eigenbasis of the symmetric matrix $\mathbf L(\xx^*)$. Furthermore, for each eigenvalue $\lambda_i(\xx^*)$, we define $\TT_i\in \real^{n\times n}$ with a permutation so that the eigenvectors associated with $\lambda_i(\xx^*)$ appear in the last columns of the matrix. As a consequence, we can define the similar matrix
$$
\mathbf{D}^{(i)}(\xx) :=\TT_i^\top \mathbf L(\xx)\TT_i.
$$
Note importantly that the matrix $\TT_i$ is defined in relation to the state $\xx^*$ and is constant for all $\xx$, so $\mathbf{D}^{(i)}$ is continuous. On the other hand, $\TT_i$ being constant does not guarantee that $\mathbf{D}^{(i)}$ will be diagonal at states other than $\xx^*$. Furthermore, by defining the matrix $\BB^{(i)}(\xx)  := \mathbf{D}^{(i)}(\xx)-\lambda_i(\xx)\mathbf I$, we can equivalently write each eigenequation with
\begin{equation}\label{eq:eig_individual}
\begin{bmatrix}
\BB^{(i)}_{aa}(\xx) & \BB^{(i)}_{ab}(\xx)\\
\BB^{(i)}_{ab}(\xx)^\top & \BB^{(i)}_{bb}(\xx) 
\end{bmatrix} \begin{bmatrix}\ww_{i,a} \\ \ww_{i,b}\end{bmatrix} =0 ,
\end{equation}
where $\ww_i$ is the vector~$\xi_i$ in the coordinate system $\TT_i$, i.e., $\xi_i=\TT_i\ww_i$. Above, we partition the matrix $\BB^{(i)}$ and the vector $\ww_i$ so that $\ww_{i,b}$ has the same dimension as the eigenspace associated with $\lambda_i$ at $\xx^*$. The next result shows that each individual eigenspace, when normalized, is already UHC. 
\begin{lemma}\label{lem:UHC_individual}\longthmtitle{UHC of individual eigenspaces}
Consider a continuous function $\map{\mathbf L}{\real^N}{\sym_n}$.
Given a state $\xx^*$ and $\delta_\ww>0$, there exists $\delta_\xx>0$ small enough such that if $\xx\in \Ball{\delta_\xx}{\xx^*}$, then for any $\ww_i$ satisfying $\BB^{(i)}(\xx)\ww_i = 0$, there exists $\ww^*$ satisfying $\BB^{(i)}(\xx^*)\ww^* = 0$ with $\|\ww_i-\ww^*\|<\delta_\ww\|\ww_i\|$.
\end{lemma}
\begin{pf}
Because $\BB^{(i)}_{aa}(\xx^*)$ is invertible, there exists $\bar \delta_{\xx}>0$ such that $\BB^{(i)}_{aa}(\xx)$ remains invertible for each $\xx \in \Ball{\bar \delta_\xx}{\xx^*}$. From~\eqref{eq:eig_individual},
$$
\ww_{i,a} = \BB^{(i)}_{aa}(\xx)^{-1}\BB^{(i)}_{ab}(\xx)\ww_{i,b}. 
$$
Because $\BB^{(i)}_{aa}(\xx)^{-1}\BB^{(i)}_{ab}(\xx)$ is continuous on $\Ball{\bar \delta_\xx}{\xx^*}$, given $\delta_\ww$, 
there exists $0<\delta_\xx<\bar \delta_\xx$ such that $\|\BB^{(i)}_{aa}(\xx)^{-1}\BB^{(i)}_{ab}(\xx)\|_F <\delta_\ww/\sqrt{2}$ for all $\xx \in \Ball{\delta_\xx}{\xx^*}$. Then,
\begin{align*}
\|\ww_{i,a}\| \leq \|\BB^{(i)}_{aa}(\xx)^{-1}\BB^{(i)}_{ab}(\xx)\|_F\|\ww_{i,b}\| < \delta_\ww\|\ww_i\|/\sqrt{2}.
\end{align*}
This also implies $\|\ww_{i,b}\|>\|\ww_i\|(1-\delta_\ww/\sqrt{2})$. Let $\ww^* = \begin{bmatrix}
0 & \ww_{i,b}^\top\|\ww_i\|/\|\ww_{i,b}\|
\end{bmatrix}^\top$ (and $\ww^*=0$ if $\|\ww_b\|=0$), then $\BB^{(i)}(\xx^*)\ww^* = 0$ because $\BB^{(i)}_{ab}(\xx^*)$ and $\BB^{(i)}_{bb}(\xx^*)$ are zero by construction. Also, we can bound the distance
\begin{align*}
\|\ww_i-\ww^*\| &= \left\Vert\begin{bmatrix}
\ww_{i,a} \\ \ww_{i,b}
\end{bmatrix}-\begin{bmatrix}
0 \\ \ww_{i,b}\|\ww_i\|/\|\ww_{i,b}\|
\end{bmatrix}\right\Vert \\
&= (\|\ww_{i,a}\|^2 + \|\ww_{i,b}\|^2(1-\|\ww_i\|/\|\ww_{i,b}\|)^2)^{1/2} \\
&= (\|\ww_{i,a}\|^2 + (\|\ww_{i,b}\|-\|\ww_i\|)^2)^{1/2}\\
&< \delta_\ww\|\ww_i\|,
\end{align*}
and the proof concludes.~\hfill~\qed
\end{pf}

From this result for individual eigenspaces, we can deduce further that any merged eigenspace is UHC. 
\begin{theorem}\longthmtitle{UHC of Merged Eigenspaces}
Consider a continuous function $\map{\mathbf L}{\real^N}{\sym_n}$. For any $\mathcal I \subseteq \until{n}$, the merged eigenspace~$\VV_\mathcal{I}$ is UHC.
\end{theorem}
\begin{pf}
Given any $\vv\in\VV_{\mathcal I}(\xx)$ in~\eqref{eq:eigspan}, we assume, without loss of generality, that if $\lambda_j (\xx)= \lambda_i(\xx)$ for some $i>j$, the associated eigenvector~$\xi_j$ is zero. This way, there is only one nonzero vector $\xi_i$ from each eigenspace. In addition, by scaling $\xi_i$, we can assume $\cc=\ones$. Using these simplifications, $\|\xi_i\|\leq 1$ because of the orthogonality of eigenspaces and the fact $\|\vv\|=1$. Thus, when we transform the coordinate frame $\ww_i = \TT_i^\top\xi_i$, we also guarantee~$\|\ww_i\|\leq 1$. This is particularly useful when we apply Lemma~\ref{lem:UHC_individual} as follows.

Consider any arbitrary $\xx^*$ at which we wish to prove UHC for $\VV_\mathcal{I}$. Lemma~\ref{lem:UHC_individual} guarantees for any given $\delta_\ww>0$ the existence of a small enough neighborhood $\Ball{\delta_\xx}{\xx^*}$ such that for every $\xx\in\Ball{\delta_\xx}{\xx^*}$, any $\ww_i$ satisfying $\BB^{(i)}(\xx)\ww_i=0$ has a corresponding $\ww^*_i\in\Ball{\delta_\ww}{\ww}$ satisfying $\BB^{(i)}(\xx^*)\ww^*_i=0$. Through coordinate transformation $\xi_i^*=\TT_i\ww_i^*$, we deduce that given the set of vectors $\{\xi_i\}_{i\in \mathcal I}$ defining $\vv$, there exists a corresponding set of vectors $\{\xi_i^*\}_{i\in \mathcal I}$ such that $\xi_i\in \Ball{\delta_\ww}{\xi_i}$ and $(\mathbf L(\xx)-\lambda_i(\xx)\mathbf I)\xi_i^* = 0$. We then define $\vv^* = (\sum \xi_i^*)/\|\sum\xi_i^*\|$, which is an element of the set $\VV_{\mathcal I}(\xx^*)$ by definition. 

We next prove that $\vv^*$ is close to $\vv$ for a small enough $\delta_\ww$. From the condition $1 = \|\vv\|= \|\sum (\xi_i^*+(\xi_i-\xi_i^*))\|$, we can bound the norm
$\|\sum\xi_i^*\| \in (1-n\delta_\ww,1+n\delta_\ww)$.  With these facts, we bound the distance
\begin{align*}
\|\vv-\vv^*\| &= \left\Vert\sum \xi_i- \frac{\sum \xi_i^*}{\|\sum\xi_i^*\|}\right\Vert \\
&\leq \|\sum (\xi_i-\xi_i^*)\|+\left\Vert\sum \left(\xi_i^*-\frac{\xi_i^*}{\|\sum\xi_i^*\|}\right)\right\Vert\\
&\leq n \delta_\ww+(1+n\delta_\ww)n\delta_\ww/(1-n\delta_\ww)\\
&= 2n \delta_\ww/(1-n\delta_\ww).
\end{align*}
Given any $\delta_\vv$, we can pick $\delta_\ww$ small enough so that  $\|\vv-\vv^*\| <\delta_\vv$, i.e., $\vv\in \Ball{\delta_\vv}{\vv^*}$. 

We have shown that given any $\delta_\vv>0$, there exists $\delta_\xx>0$ such that any $\vv\in\VV_{\mathcal I}(\xx)$, for $\xx\in\Ball{\delta_\xx}{\xx^*}$, has a corresponding $\vv^*\in\VV_{\mathcal I}(\xx^*)$ such that $\vv\in\Ball{\delta_\vv}{\vv^*}$. In other words, $\VV_{\mathcal I}(\xx)$, for $\xx\in\Ball{\delta_\xx}{\xx^*}$, is a subset of a $\delta_\vv$ neighborhood of $\VV_{\mathcal I}(\xx^*)$, which is precisely the definition of UHC, concluding the proof.~\hfill~\qed
\end{pf}

\subsection{Lower Hemicontinuity of Merged Eigenspaces}
Unlike the case of UHC, individual normalized eigenspaces are not LHC everywhere. Therefore, we proceed directly to the analysis of the merged eigenspaces. We define, for an index set $\mathcal I \subseteq \until n$, an orthonormal eigenbasis matrix~$\TT_{\mathcal I}\in \real^{n\times n}$, with the eigenvectors associated with $\lambda_i(\xx^*)$ for $i\in \mathcal I$ showing up in the last columns of the matrix. Then, we define the matrix
$$
\mathbf{D}^{\mathcal I}(\xx) = \TT_{\mathcal I}^\top \mathbf L(\xx)\TT_{\mathcal I}.
$$
The next result establishes the LHC property of the merged eigenspaces.

\begin{theorem}\longthmtitle{LHC of Merged Eigenspaces} Consider a continuous function $\map{\mathbf L}{\real^N}{\sym_n}$. For any $\mathcal I \subseteq \until n$, the merged eigenspace is LHC at $\xx$ where $\lambda_i(\xx)\neq \lambda_j(\xx)$ for all $i\in \mathcal I$ and $j\not \in \mathcal I$, i.e., where none of the eigenvalues considered in the span is equal to any of the remaining eigenvalue.
\end{theorem}
\begin{pf}
Consider the change of coordinate frame $\xi_i = \TT_{\mathcal I} \ww_i$, for each $i\in\until{n}$. The merged eigenspace given by \eqref{eq:eigspan} can be rewritten as
\begin{multline*}
  \VV_{\mathcal I}(\xx) = \setdefB{\vv\in \real^n}{
  (\mathbf{D}^{\mathcal I}(\xx)-\lambda_i(\xx)\mathbf I)\ww_i = 0,~\forall i\in \mathcal I, \\
  \vv = \TT_{\mathcal I}\WW\cc,~\cc\in\real^{\vert \mathcal I\vert},~\|\vv\|=1 } ,
\end{multline*}
where $\WW\in \real^{n\times \vert\mathcal I\vert}$ is a matrix constructed by stacking $\ww_i$ together. By construction, given an element $\vv^*\in \VV_{\mathcal I}(\xx^*)$, it must take the form $\vv^*= \TT_{\mathcal I}\left[\begin{smallmatrix}
0 \\ \psi
\end{smallmatrix}\right]$ for some $\psi\in \real^{\vert \mathcal I \vert}$. 

Consider $\xx^*$ at which we wish to prove LHC for~$\VV_\mathcal{I}$. We next show the existence of $\vv \in \VV_{\mathcal I}(\xx)$ close enough to $\vv^*$ for all $\xx$ close enough to $\xx^*$. First, we partition the eigenequations,
$$
\left(\begin{bmatrix}
\mathbf{D}_{aa}^{\mathcal I}(\xx) & \mathbf{D}_{ab}^{\mathcal I}(\xx) \\
\mathbf{D}_{ab}^{\mathcal I}(\xx)^\top & \mathbf{D}_{bb}^{\mathcal I}(\xx)
\end{bmatrix}-\lambda_i(\xx) \mathbf I\right)\begin{bmatrix}
\ww_{i,a} \\ \ww_{i,b}
\end{bmatrix} = 0,
$$
so that $\ww_{i,b}$ has the dimension of $\vert \mathcal I\vert$. The matrix $\mathbf{D}_{aa}^{\mathcal I}(\xx^*)$ is a diagonal matrix of eigenvalues $\lambda_j(\xx^*)$ for $j\not \in \mathcal{I}$. Because $\lambda_j(\xx^*)\neq\lambda_i(\xx^*)$ for any $i\in \mathcal{I}$ and $j\not \in \mathcal{I}$, the matrix $\mathbf{D}_{aa}^{\mathcal I}(\xx)- \lambda_i(\xx) \mathbf I$ is invertible at $\xx=\xx^*$. Then due to continuity of the matrix, there exists $\bar \delta_\xx$ such that it remains invertible for $\xx \in \Ball{\bar\delta_\xx}{\xx^*}$, and we can find the following relationship,
$$
\ww_{i,a} = (\mathbf{D}^{\mathcal I}_{aa}(\xx)-\lambda_{i}(\xx)\mathbf I)^{-1}\mathbf{D}^{\mathcal I}_{ab}(\xx)\ww_{i,b}.
$$
Due to continuity of the matrix $\mathbf{D}^\mathcal{I}$ and the fact that $\mathbf{D}^{\mathcal I}_{ab}(\xx)$ is a zero matrix at $\xx=\xx^*$, we can further find that given $\delta_\ww$, there exists $0<\delta_\xx\leq \bar \delta_\xx$ such that $\|\ww_{i,a}\| \leq \delta_\ww \|\ww_{i,b}\|$ for all $\xx \in \Ball{\delta_\xx}{\xx^*}$. With this property, we construct $\vv\in\VV_{\mathcal I}(\xx)$ with the following procedure.

We begin by selecting the set of eigenvectors $\{\ww_i\}_{i\in\mathcal{I}}$ to be orthonormal to one another. This set of eigenvectors must exist because $\mathbf{D}^{\mathcal I}(\xx)$ is symmetric. With this choice, we can show that when we partition the matrix $\WW =  \left[\begin{smallmatrix}
\WW_a \\ \WW_b
\end{smallmatrix}\right]$, $\WW_b$ is an invertible matrix for $\xx \in \Ball{\delta_\xx}{\xx^*}$. We prove this statement by contradiction. Assume that $\WW_b$ is not full rank, then there exists a vector $0\neq \cc\in \real^{\vert\mathcal I\vert}$ such that $\WW_b\cc = 0$. In addition, $\WW^\top \WW = \mathbf I$ because $\ww_i$ are orthogonal to each other. Thus,
$$
\|\cc\| =\|\WW^\top \WW\cc\| = \|\WW_a^\top \WW_a \cc\| \leq \delta_\ww^2 \vert \mathcal I \vert^2 \|\cc\|, 
$$
which is a contradiction for small $\delta_\ww$. Since $\WW_b$ is invertible, we can define the vector
$$\bar \vv = \TT_{\mathcal I}\WW\WW_b^{-1}\psi = \TT_{\mathcal I}\begin{bmatrix}
\WW_a\WW_b^{-1}\psi \\ \psi
\end{bmatrix},$$
which we use to construct $\vv\in\VV_{\mathcal I}(\xx)$
Before doing so, we upper bound $\WW_a\WW_b^{-1}\psi$. Note that
\begin{align*}
\|\WW_a\WW_b^{-1}\psi\|&=\|\WW_a(W^\top_b\WW_b)^{-1}\WW_b^\top\psi\| \\
&\leq \|\WW_a\|\|(\WW^\top_b\WW_b)^{-1}\|\|\WW_b\|.
\end{align*}
Here, we can bound $\|\WW_b\|\leq \vert\mathcal I\vert$ due to normality of each $\ww_i$. Also from the earlier fact $\|\ww_{i,a}\| \leq \delta_\ww \|\ww_{i,b}\|\leq \delta_\ww$, we bound $\|\WW_a\|\leq \delta_\ww\vert\mathcal I\vert$. As for the $\|(\WW^\top_b\WW_b)^{-1}\|$, we investigate the smallest eigenvalue of $(\WW^\top_b\WW_b)$. Due to orthonormality,
$$
\ww^\top_{i,b}\ww_{j,b} =
\begin{cases}
 - \ww^\top_{i,a}\ww_{j,a} & i\neq j ,\\
 1- \ww^\top_{i,a}\ww_{j,a} & i=j .
\end{cases}
$$
Combined with the fact $\|\ww_{i,a}\|\leq \delta_\ww$, we upper bound the off-diagonal entries of $\WW^\top_b\WW_b$ with $\delta_\ww^2$, and we lower bound the diagonal entries with  $1-\delta_\ww^2$. Using the Gershgorin circle theorem~\citep[Thm. 1.3]{FB-JC-SM:09}, the smallest eigenvalue of $\WW^\top_b\WW_b$ is lower bounded by $1-\vert \mathcal I \vert\delta_\ww^2$. Using these bounds, we find
$$
\|\WW_a\WW_b^{-1}\psi\| \leq \frac{\delta_\ww\vert \mathcal I \vert}{1-\delta_\ww^2} := \delta_\vv.
$$
Note here that smaller $\delta_\vv$ corresponds to small $\delta_\ww<1$.

Finally, we select $\cc=\ww_b^{-1}\psi/\|\bar \vv\|$ to construct $\vv = \bar \vv / \|\bar \vv\|$, which is an element of $\VV_{\mathcal I}(\xx)$. Let $\theta$ be the angle between the unit vectors $\vv$ and $\vv^*$. Then, we bound
\begin{align*}
    \|\vv-\vv^*\| \leq \theta \leq \tan \theta = \frac{\|\WW_a\WW_b^{-1}\psi\|}{\|\psi\|} \leq \delta_\vv.
\end{align*}
Thus, we have proven that given any $\delta_\vv>0$, there exists $\delta_\xx>0$ such that if $\xx \in \Ball{\delta_\xx}{\xx^*}$, then there exists $\vv\in\VV_{\mathcal I}(\xx)$ where $\vv\in \Ball{\delta_\vv}{\vv^*}$. This is sufficient to prove that given any sequence $\{\xx^k\}_{k\in \naturals}$ converging to $\xx^*$, there exists a sequence $\{\vv^k\}_{k\in \naturals}$, with $\vv^k \in \VV_{\mathcal I}(\xx^k)$, converging to $\vv^*$, concluding the proof.~\hfill~\qed
\end{pf}


{
  \small
  \begin{thebibliography}{42}
\providecommand{\natexlab}[1]{#1}
\providecommand{\url}[1]{\texttt{#1}}
\expandafter\ifx\csname urlstyle\endcsname\relax
  \providecommand{\doi}[1]{doi: #1}\else
  \providecommand{\doi}{doi: \begingroup \urlstyle{rm}\Url}\fi

\bibitem[Aliprantis and Border(1999)]{CDA-KCB:99}
C.~D. Aliprantis and K.~C. Border.
\newblock \emph{Infinite Dimensional Analysis: A Hitchhiker's Guide}.
\newblock Studies in Economic Theory. Springer, New York, 1999.
\newblock ISBN 9783540658542.

\bibitem[Ames et~al.(2017)Ames, Xu, Grizzle, and Tabuada]{ADA-XX-JWG-PT:17}
A.~D. Ames, X.~Xu, J.~W. Grizzle, and P.~Tabuada.
\newblock Control barrier function based quadratic programs for safety critical
  systems.
\newblock \emph{IEEE Transactions on Automatic Control}, 62\penalty0
  (8):\penalty0 3861--3876, 2017.

\bibitem[Ames et~al.(2019)Ames, Coogan, Egerstedt, Notomista, Sreenath, and
  Tabuada]{ADA-SC-ME-GN-KS-PT:19}
A.~D. Ames, S~Coogan, M.~Egerstedt, G.~Notomista, K.~Sreenath, and P.~Tabuada.
\newblock Control barrier functions: Theory and applications.
\newblock In \emph{{E}uropean {C}ontrol {C}onference}, pages 3420--3431,
  Naples, Italy, June 2019.

\bibitem[Artstein(1983)]{ZA:83}
Z.~Artstein.
\newblock Stabilization with relaxed controls.
\newblock \emph{Nonlinear Analysis}, 7\penalty0 (11):\penalty0 1163--1173,
  1983.

\bibitem[Blanchini and Miani(2007)]{FB-SM:07}
F.~Blanchini and S.~Miani.
\newblock \emph{Set-Theoretic Methods in Control}.
\newblock Birkh{\"a}user, Boston, MA, 2007.
\newblock ISBN 9780817646066.

\bibitem[Border(1985)]{KCB:85}
K.~C. Border.
\newblock \emph{Fixed Point Theorems with Applications to Economics and Game
  Theory}.
\newblock Cambridge University Press, Cambridge, UK, 1985.
\newblock ISBN 9780521388085.

\bibitem[Boyd(2006)]{SB:06}
S.~Boyd.
\newblock Convex optimization of graph {L}aplacian eigenvalues.
\newblock In \emph{Int. Congress of Mathematicians}, pages 1311--1319, Madrid,
  Spain, August 2006.

\bibitem[Boyd and Vandenberghe(2009)]{SB-LV:09}
S.~Boyd and L.~Vandenberghe.
\newblock \emph{Convex Optimization}.
\newblock Cambridge University Press, Cambridge, UK, 2009.
\newblock ISBN 0521833787.

\bibitem[Bullo et~al.(2009)Bullo, Cort{\'e}s, and Martinez]{FB-JC-SM:09}
F.~Bullo, J.~Cort{\'e}s, and S.~Martinez.
\newblock \emph{Distributed Control of Robotic Networks}.
\newblock Applied Mathematics Series. Princeton University Press, Princeton,
  NJ, 2009.
\newblock ISBN 978-0-691-14195-4.

\bibitem[Capelli and Sabattini(2020)]{BC-LS:20}
B.~Capelli and L.~Sabattini.
\newblock Connectivity maintenance: Global and optimized approach through
  control barrier functions.
\newblock In \emph{{IEEE} Int. Conf.\ on Robotics and Automation}, pages
  5590--5596, Paris, France, May 2020.

\bibitem[Cherukuri et~al.(2017)Cherukuri, Gharesifard, and
  Cort{\'e}s]{AC-BG-JC:17}
A.~Cherukuri, B.~Gharesifard, and J.~Cort{\'e}s.
\newblock Saddle-point dynamics: conditions for asymptotic stability of saddle
  points.
\newblock \emph{SIAM Journal on Control and Optimization}, 55\penalty0
  (1):\penalty0 486--511, 2017.

\bibitem[Clarke(1983)]{FHC:83}
F.~H. Clarke.
\newblock \emph{Optimization and Nonsmooth Analysis}.
\newblock Canadian Mathematical Society Series of Monographs and Advanced
  Texts. Wiley, New York, 1983.
\newblock ISBN 047187504X.

\bibitem[Cort\'es and Egerstedt(2017)]{JC-ME:17-jcmsi}
J.~Cort\'es and M.~Egerstedt.
\newblock Coordinated control of multi-robot systems: A survey.
\newblock \emph{SICE Journal of Control, Measurement, and System Integration},
  10\penalty0 (6):\penalty0 495--503, 2017.

\bibitem[de~Gennaro and Jadbabaie(2006)]{MDG-AJ:06}
M.~C. de~Gennaro and A.~Jadbabaie.
\newblock Decentralized control of connectivity for multi-agent systems.
\newblock In \emph{{IEEE} Conf.\ on Decision and Control}, pages 3628--3633,
  San Diego, CA, December 2006.

\bibitem[Egerstedt et~al.(2018)Egerstedt, Pauli, Notomista, and
  Hutchinson]{ME-JNP-GN-SH:18}
M.~Egerstedt, J.~N. Pauli, G.~Notomista, and S.~Hutchinson.
\newblock Robot ecology: Constraint-based control design for long duration
  autonomy.
\newblock \emph{Annual Reviews in Control}, 46:\penalty0 1--7, 2018.

\bibitem[Fiedler(1973)]{MF:73}
M.~Fiedler.
\newblock Algebraic connectivity of graphs.
\newblock \emph{Czechoslovak Mathematical Journal}, 23\penalty0 (2):\penalty0
  298--305, 1973.

\bibitem[Freeman and Kototovic(1996)]{RAF-PVK:96}
R.~A. Freeman and P.~V. Kototovic.
\newblock \emph{Robust Nonlinear Control Design: State-space and Lyapunov
  Techniques}.
\newblock Birkh{\"a}user, Boston, MA, 1996.
\newblock ISBN 0-8176-3930-6.

\bibitem[Gasparri et~al.(2017)Gasparri, Sabattini, and Ulivi]{AG-LS-GU:17}
A.~Gasparri, L.~Sabattini, and G.~Ulivi.
\newblock Bounded control law for global connectivity maintenance in
  cooperative multi-robot systems.
\newblock \emph{IEEE Transactions on Robotics}, 33\penalty0 (3):\penalty0
  700--717, 2017.

\bibitem[Glotfelter et~al.(2017)Glotfelter, Cort\'es, and
  Egerstedt]{PG-JC-ME:17}
P.~Glotfelter, J.~Cort\'es, and M.~Egerstedt.
\newblock Nonsmooth barrier functions with applications to multi-robot systems.
\newblock \emph{IEEE Control Systems Letters}, 1\penalty0 (2):\penalty0
  310--315, 2017.

\bibitem[Godsil and Royle(2001)]{CDG-GFR:01}
C.~D. Godsil and G.~F. Royle.
\newblock \emph{Algebraic Graph Theory}, volume 207 of \emph{Graduate Texts in
  Mathematics}.
\newblock Springer, 2001.
\newblock ISBN 0387952411.

\bibitem[Hale(1969)]{JKH:69}
J.~K. Hale.
\newblock \emph{Ordinary Differential Equations}.
\newblock Wiley, New York, 1969.
\newblock ISBN 9780471340904.

\bibitem[Ji and Egerstedt(2007)]{MJ-ME:07}
M.~Ji and M.~Egerstedt.
\newblock Distributed coordination control of multiagent systems while
  preserving connectedness.
\newblock \emph{IEEE Transactions on Robotics}, 23\penalty0 (4):\penalty0
  693--703, 2007.

\bibitem[Kato(1976)]{TK:76}
T.~Kato.
\newblock \emph{Perturbation Theory for Linear Operators}.
\newblock Grundlehren der mathematischen Wissenschaften: a series of
  comprehensive studies in mathematics. Springer, Berlin, 1976.

\bibitem[Kim and Mesbahi(2006)]{YK-MM:06}
Y.~Kim and M.~Mesbahi.
\newblock On maximizing the second smallest eigenvalue of a state-dependent
  graph {L}aplacian.
\newblock \emph{IEEE Transactions on Automatic Control}, 51\penalty0
  (1):\penalty0 116--120, 2006.

\bibitem[Lechicki and Spakowski(1985)]{AL-AS:85}
A.~Lechicki and A.~Spakowski.
\newblock A note on intersection of lower semicontinuous multifunctions.
\newblock \emph{Proc. of the American Mathematical Society}, 95\penalty0
  (1):\penalty0 119--122, 1985.

\bibitem[Lewis(1996)]{ASL:96}
A.~S. Lewis.
\newblock Group invariance and convex matrix analysis.
\newblock \emph{SIAM Journal on Matrix Analysis and Applications}, 17\penalty0
  (4):\penalty0 927--949, 1996.

\bibitem[Mesbahi and Egerstedt(2010)]{MM-ME:10}
M.~Mesbahi and M.~Egerstedt.
\newblock \emph{Graph Theoretic Methods in Multiagent Networks}.
\newblock Applied Mathematics Series. Princeton University Press, 2010.

\bibitem[Morris et~al.(2015)Morris, Powell, and Ames]{BJM-MJP-ADA:15}
B.~J. Morris, M.~J. Powell, and A.~D. Ames.
\newblock Continuity and smoothness properties of nonlinear optimization-based
  feedback controllers.
\newblock In \emph{{IEEE} Conf.\ on Decision and Control}, pages 151--158,
  Osaka, Japan, Dec 2015.

\bibitem[\"Ogren et~al.(2006)\"Ogren, Backlund, Harryson, Kristensson, and
  Stensson]{PO-AB-TH-LK-PS:06}
P.~\"Ogren, A.~Backlund, T.~Harryson, L.~Kristensson, and P.~Stensson.
\newblock Autonomous {UCAV} strike missions using behavior control {L}yapunov
  functions.
\newblock In \emph{AIAA Guidance, Navigation, and Control Conference and
  Exhibit}, page 6197, August 2006.

\bibitem[Ong et~al.(2021)Ong, Capelli, Sabattini, and
  Cort\'{e}s]{PO-BC-LS-JC:21}
P.~Ong, B.~Capelli, L.~Sabattini, and J.~Cort\'{e}s.
\newblock Network connectivity maintenance via nonsmooth control barrier
  functions.
\newblock In \emph{{IEEE} Conf.\ on Decision and Control}, pages 4780--4785,
  Austin, TX, December 2021.

\bibitem[Oriolo et~al.(2002)Oriolo, Luca, and Vendittelli]{GO-ADL-MV:02}
G.~Oriolo, A.~De Luca, and M.~Vendittelli.
\newblock {WMR} control via dynamic feedback linearization: design,
  implementation, and experimental validation.
\newblock \emph{IEEE Transactions on Control Systems Technology}, 10\penalty0
  (6):\penalty0 835--852, 2002.

\bibitem[Prajna and Jadbabaie(2004)]{SP-AJ:04}
S.~Prajna and A.~Jadbabaie.
\newblock Safety verification of hybrid systems using barrier certificates.
\newblock In \emph{Hybrid Systems: Computation and Control}, pages 477--492,
  Philadelphia, PA, March 2004.

\bibitem[Rockafellar(1970)]{RTR:70}
R.~T. Rockafellar.
\newblock \emph{Convex Analysis}.
\newblock Princeton University Press, Princeton, NJ, 1970.
\newblock ISBN 9780691015866.

\bibitem[Sabattini et~al.(2013)Sabattini, Chopra, and Secchi]{LS-NC-CS:13}
L.~Sabattini, N.~Chopra, and C.~Secchi.
\newblock Decentralized connectivity maintenance for cooperative control of
  mobile robotic systems.
\newblock \emph{International Journal of Robotics Research}, 32\penalty0
  (12):\penalty0 1411--1423, 2013.

\bibitem[Schuresko and Cort{\'e}s(2009)]{MDS-JC:09}
M.~D. Schuresko and J.~Cort{\'e}s.
\newblock Distributed motion constraints for algebraic connectivity of robotic
  networks.
\newblock \emph{Journal of Intelligent and Robotic Systems}, 56\penalty0
  (1):\penalty0 99--126, 2009.

\bibitem[Schuresko and Cort{\'e}s(2012)]{MDS-JC:09-sicon}
M.~D. Schuresko and J.~Cort{\'e}s.
\newblock Distributed tree rearrangements for reachability and robust
  connectivity.
\newblock \emph{SIAM Journal on Control and Optimization}, 50\penalty0
  (5):\penalty0 2588--2620, 2012.

\bibitem[Shevitz and Paden(1994)]{DS-BP:94}
D.~Shevitz and B.~Paden.
\newblock Lyapunov stability theory of nonsmooth systems.
\newblock \emph{IEEE Transactions on Automatic Control}, 39\penalty0
  (9):\penalty0 1910--1914, 1994.

\bibitem[Still(2018)]{GS:18}
G.~Still.
\newblock Lectures on parametric optimization: An introduction.
\newblock \emph{Optimization Online}, 2018.

\bibitem[Warner(1989)]{FWW:89}
F.~W. Warner.
\newblock \emph{Foundations of Differentiable Manifolds and Lie Groups}.
\newblock Number~94 in Graduate Texts in Mathematics. Springer, New York, 2
  edition, 1989.
\newblock ISBN 9780387908946.

\bibitem[Wieland and Allg{\"o}wer(2007)]{PW-FA:07}
P.~Wieland and F.~Allg{\"o}wer.
\newblock Constructive safety using control barrier functions.
\newblock \emph{IFAC Proceedings Volumes}, 40\penalty0 (12):\penalty0 462--467,
  2007.

\bibitem[Zavlanos and Pappas(2005)]{MMZ-GJP:05}
M.~M. Zavlanos and G.~J. Pappas.
\newblock Controlling connectivity of dynamic graphs.
\newblock In \emph{{IEEE} Conf.\ on Decision and Control and European Control
  Conference}, pages 6388--6393, Seville, Spain, December 2005.

\bibitem[Zavlanos and Pappas(2015)]{MZZ-GJP:14-sv}
M.~M. Zavlanos and G.~J. Pappas.
\newblock Connectivity of dynamic graphs.
\newblock In J.~Baillieul and T.~Samad, editors, \emph{Encyclopedia of Systems
  and Control}. Springer-Verlag, 2015.

\end{thebibliography}

}
\end{document}